\documentclass[10pt,twoside,leqno]{amsart}
\usepackage{amssymb,amsbsy,amsmath,amsfonts,amssymb,amscd,times,graphics,color,wasysym}
\newcommand{\C}                 {\mathbb{C}}

\newcommand{\K}                 {\mathbb{K}}
\newcommand{\R}                 {\mathbb{R}}
\newcommand{\N}                 {\mathbb{N}}

\newtheorem{theorem}{Theorem}
\newtheorem{lemma}{Lemma}

\begin{document}

%\large

% Dans \equation \small : 
% Ctrl eo : \underline{ }_{ \octagon \! \! \! \! \! \! \tiny{ }} 
% Ctrl ef : \underline{ }_{ \fbox{\tiny }}
% Ctrl ee : \frac{ 1}{ 2} \,
% Ctrl er : \delta_{ l_1}^j \,

%\def\theequation{*}\begin{equation}
%\underline{ 
%-\frac{ 1}{ 6} \, L_{l_1, l_1}^{l_1} \, \Theta_x^0
%}_{ 
%\octagon \! \! \! \! \! \! \tiny{ a}} 
%\end{equation}

%\def\theequation{**}\begin{equation}
%\underline{
%\frac{ 1}{ 2} \, H_{ l_1}^j \, \Theta^0 
%}_{ 
%\fbox{\tiny 1}}
%\end{equation}

\title[
Differential characterization of $Y_{ X^{j_1} X^{ j_2} } = 0$, 
$1\leq j_1, j_2 \leq n\geq 2$
]{
Explicit differential characterization of
\\ 
PDE systems pointwise equivalent to $Y_{ X^{j_1} X^{ j_2} } = 0$
\\
$1\leq j_1, j_2 \leq n\geq 2$.}

\author{Jo\"el Merker}

\address{
CNRS, Universit\'e de Provence, LATP, UMR 6632, CMI, 
39 rue Joliot-Curie, F-13453 Marseille Cedex 13, France}

\email{merker@cmi.univ-mrs.fr} 


\date{\number\year-\number\month-\number\day}

\begin{abstract} 
In~\cite{ lie1883}, as an early result, Sophus Lie established 
that a second order ordinary differential 
equation $y_{xx }=F( x, y, y_x)$ is
equivalent, through an invertible point transformation $(x, y)\mapsto (X
(x,y), Y (x,y))$, to the free particle equation $Y_{ XX } =0$ 
{\it if and only if}\, the right member $F$ is a degree three
polynomial in $y_x$, namely there exist four functions $G$, $H$, $L$
and $M$ of $(x,y)$ such that $F$ can be written as
$$
F(x,y,y_x)=G(x,y)+y_x\cdot H(x,y)+
(y_x)^2\cdot L(x,y)+(y_x)^3\cdot M(x,y),
$$
and furthermore, the four functions $G$, $H$, $L$ and $M$ satisfy two
\underline{second order} partial differential equations:
$$
\aligned
0 
&
= 
-
2G_{yy}
+
\frac{4}{3}H_{xy}
-
\frac{2}{3}L_{xx}
+
2(GL)y
-
2G_xM
-
4GM_x
+
\frac{2}{3}HL_x
-
\frac{4}{3}
HH_y,
\\
0
&
=
-
\frac{2}{3}H_{yy}
+
\frac{4}{3}L_{xy}
-
2M_{xx}
+
2GM_y
+
4G_yM
-
2(HM)_x
-
\frac{2}{3}H_yL
+
\frac{4}{3}LL_x.
\endaligned
$$
In~\cite{ m2004a}, this theorem was generalized to systems of Newtonian
particles with $m\geq 2$ degree of freedom, {\it i.e.} with one
independent variable $x$ and $m\geq 2$ dependent variables $(y^1, y^2,
\dots, y^m)$. In this paper, we generalize S.~Lie's theorem to the
case of several independent variables $(x^1, x^2 \dots, x^n)$, $n\geq
2$, and one dependent variable $y$. Strikingly, as in~\cite{ m2004a},
the (complicated) differential system which corresponds to the above
two second order partial differential equations is of \underline{first
order}. By means of computer programming, this phenomenon was
discovered in~\cite{ bn2002}, \cite{ n2003} in the case $n=2$.
In~\cite{ bi2003}, \cite{ m2003}, the general case $n\geq 2$ was
handled.
\end{abstract}

\maketitle

\begin{center}
\begin{minipage}[t]{11cm}
\baselineskip =0.35cm
{\scriptsize

\centerline{\bf Table of contents}

\smallskip

{\bf 1.~Introduction \dotfill 1.}

{\bf 2.~Completely integrable systems of second order ordinary
differential equations \dotfill 9.}

{\bf 3.~First and second auxiliary systems \dotfill 16.} 

}\end{minipage}
\end{center}

\section*{\S1.~Introduction}

This paper, a direct continuation of~\cite{ m2004a}, provides a
summarized proof of the following statement, labelled as Theorem~1.23
in the Introduction of~\cite{ m2004a}. All our functions are assumed to
be analytic.

\def\thetheorem{1.1}\begin{theorem}
{\rm ($n=2$: \cite{bn2002}, \cite{ n2003}; 
$n\geq 2$: \cite{ bi2003}, \cite{ m2003})} \ Let $\K = \R$ or $\C$,
let $n\in\N$, \underline{\rm suppose $n\geq 2$} and consider a system
of completely integrable partial differential equations in $n$
independent variables $x= (x^1, \dots, x^n) \in \K^n$ and in one
dependent variable $y \in \K$ of the form{\rm :}
\def\theequation{1.2}\begin{equation}
y_{x^{j_1}x^{j_2}}(x) 
=
F^{j_1, j_2} 
\left(
x, y(x), y_{x^1}(x),\dots,y_{x^n}(x)
\right), 
\ \ \ \ \ \ \ \ 
j_1,j_2= 1,\dots n, 
\end{equation}
where $F^{j_1, j_2} = F^{ j_2, j_1}$. Under a local change of
coordinates $(x, y) \mapsto (X, Y) = (X(x, y), Y(x, y))$, this
system~\thetag{1.2} is equivalent to the {\rm simplest} system
$Y_{X^{j_1} X^{j_2}}=0$, $j_1, j_2 =1, \dots, n$, {\rm if and only if}
there exist {\rm arbitrary} functions $G_{j_1, j_2}$, $H_{ j_1, j_2}^{
k_1}$, $L_{ j_1}^{ k_1}$ and $M^{k_1}$ of the variables $(x, y)$, for
$1\leq j_1,j_2,k_1 \leq n$, satisfying the
two symmetry conditions $G_{j_1, j_2}= G_{j_2, j_1}$ and $H_{j_1,
j_2}^{ k_1}= H_{j_2,j_1}^{k_1}$, such that the equation~\thetag{1.2}
is of the specific cubic polynomial form{\rm :}
\def\theequation{1.3}\begin{equation}
y_{x^{j_1}x^{j_2}}
=
G_{j_1,j_2}+\sum_{k_1=1}^n\, 
y_{x^{k_1}}\left(
H_{j_1,j_2}^{k_1} +
\frac{ 1}{2}\, y_{x^{j_1}} \, L_{j_2}^{k_1}+
\frac{ 1}{2}\, y_{x^{j_2}}\, L_{j_1}^{k_1} +
y_{x^{j_1}}
y_{x^{j_2}}\, M^{k_1}
\right), 
\end{equation}
for $j_1, j_2 = 1, \dots, n$.
\end{theorem}

We refer the reader to~\cite{ m2004a} for an extensive introduction
and for a more substantial bibliography. Applying \'E.~Cartan's
equivalence algorithm ({\it see}~\cite{ g1989} and \cite{ ol1995} for
modern expositions), the general case $n\geq 2$ of the above theorem
has also been established in~\cite{ ha1937}, where the construction of
a projective connection associated to a second order ordinary
differential equation achieved in~\cite{ ca1924} was extended to
several variables. Theorem~1.1 was re-discovered in~\cite{ bn2002}, in
\cite{ n2003}, in~\cite{ bi2003} and in~\cite{ m2003}, thanks to
partial parametric computations of the Hachtroudi-Chern tensors in
$n\geq 2$ variables, which were described in a non-parametric way
in~\cite{ ch1975} ({\it see}~\S1.8 below).

It may seem quite paradoxical and counter-intuitive (or even false?) that
{\it every} system~\thetag{ 1.3}, for {\it arbitrary} choices of
functions $G_{j_1, j_2}$, $H_{j_1, j_2}^{k_1}$, $L_{j_1}^{k_1}$ and
$M^{k_1}$, is {\it automatically} equivalent to $Y_{X^{j_1} X^{j_2}}=
0$. However, a strong hidden assumption holds: that of {\sl complete
integrability}. Shortly, this crucial condition amounts to say
that
\def\theequation{1.4}\begin{equation}
D_{x^{j_3}}\left(
F^{j_1,j_2}
\right)=D_{x^{j_2}}\left(
F^{j_1,j_3}
\right),
\end{equation}
for all $j_1, j_2, j_3 = 1, \dots, n$, where, for $j =1, \dots, n$, the
$D_{ x^j}$ are the {\sl total differentiation operators} defined by
\def\theequation{1.5}\begin{equation}
D_{x^j}:=
\frac{ \partial }{\partial x^j}+
y_{x^j} \, 
\frac{ \partial }{\partial y}+
\sum_{l=1}^n\,
F^{j,l}\, 
\frac{ \partial }{\partial y_{x^l}}.
\end{equation}
These conditions are non-void precisely when $n\geq 2$. More
concretely, writing out~\thetag{ 1.4} when the
$F^{j_1, j_2}$ are of the specific cubic polynomial form~\thetag{
1.3}, after some nontrivial manual computation, we obtain the
complicated cubic differential polynomial in the
variables $y_{x^k}$ labelled as equation~\thetag{ 1.25} 
in~\cite{ m2004a}. Equating to zero all the
coefficients of this cubic polynomial, we obtain four familes (I'),
(II'), (III') and (IV') of \underline{first order} partial
differential equations satisfied by $G_{j_1, j_2}$, $H_{j_1,
j_2}^{k_1}$, $L_{j_1}^{k_1}$ and $M^{k_1}$:
\def\theequation{I'}\begin{equation}
\left\{
0 =
G_{j_1,j_2, x^{j_3}}- G_{j_1,j_3, x^{j_2}}
+
\sum_{k_1=1}^n\, H_{j_1,j_2}^{k_1}\, G_{k_1, j_3}
-
\sum_{k_1=1}^n\, H_{j_1, j_3}^{k_1}\, G_{k_1, j_2}.
\right.
\end{equation}
\def\theequation{II'}\begin{equation}
\left\{
\aligned
0
&
=
\delta_{j_3}^{k_1}\, G_{j_1,j_2, y}
-
\delta_{j_2}^{k_1}\, G_{j_1,j_3, y}
+
H_{j_1, j_2, x^{j_3}}^{k_1}
- 
H_{j_1, j_3, x^{j_2}}^{k_1}
+ 
\\
& \
\ \ \ \ \
+
\frac{ 1}{2}\,G_{j_1, j_3}\, L_{j_2}^{k_1}
- 
\frac{ 1}{2}\,G_{j_1, j_2}\, L_{j_3}^{k_1}
+ \\
& \
\ \ \ \ \
+
\frac{1}{2}\,\delta_{j_1}^{k_1}\,\sum_{k_2=1}^n\,
G_{k_2,j_3}\,L_{j_2}^{k_2}
-
\frac{ 1}{2}\,\delta_{j_1}^{k_1}\,\sum_{k_2=1}^n\, 
G_{k_2,j_2}\, L_{j_3}^{k_2}
+ \\
& \
\ \ \ \ \
+
\frac{ 1}{2}\,\delta_{j_2}^{k_1}\,\sum_{k_2=1}^n\, 
G_{k_2, j_3}\, L_{j_1}^{k_2}
-
\frac{ 1}{2}\,\delta_{j_3}^{k_1}\,\sum_{k_2=1}^n\, 
G_{k_2,j_2}\,L_{j_1}^{k_2}
+ \\
& \
\ \ \ \ \
+
\sum_{k_2=1}^n\,H_{k_2,j_3}^{k_1}\,H_{j_1, j_2}^{k_2}
- 
\sum_{k_2=1}^n\,H_{k_2,j_2}^{k_1}\,H_{j_1,j_3}^{k_2}.
\endaligned\right.
\end{equation}
\def\theequation{III'}\begin{equation}
\left\{
\aligned
0
&
=
\sum_{\sigma\in\mathfrak{S}_2}\left(
\delta_{j_3}^{k_{2)}}\,H_{j_1,j_2,y}^{k_{\sigma(1)}}
-
\delta_{j_2}^{k_{\sigma(2)}}\,H_{j_1,j_3,y}^{k_{\sigma(1)}}
+
\right. \\
& \
\ \ \ \ \ \ \ \ \ \ \ \ \ \ \
\left.
+
\frac{1}{2}\,\delta_{j_2}^{k_{\sigma(2)}}\,L_{j_1,x^{j_3}}^{k_{\sigma(1)}}
-
\frac{1}{2}\,\delta_{j_3}^{k_{\sigma(2)}}\,L_{j_1,x^{j_2}}^{k_{\sigma(1)}}
+
\right. \\
& \
\ \ \ \ \ \ \ \ \ \ \ \ \ \ \
\left.
+
\frac{1}{2}\,\delta_{j_1}^{k_{\sigma(2)}}\,L_{j_2,x^{j_3}}^{k_{\sigma(1)}}
-
\frac{1}{2}\,\delta_{j_1}^{k_{\sigma(2)}}\,L_{j_3,x^{j_2}}^{k_{\sigma(1)}}
+
\right. \\
& \
\ \ \ \ \ \ \ \ \ \ \ \ \ \ \
\left.
+
\delta_{j_2}^{k_{\sigma(2)}}\,G_{j_1,j_3}\,M^{k_{\sigma(1)}}
-
\delta_{j_3}^{k_{\sigma(2)}}\,G_{j_1,j_2}\,M^{k_{\sigma(1)}}
+
\right.
\\
& \ 
\ \ \ \ \ \ \ \ \ \ \ \ \ \ \
\left.
+
\delta_{j_1,\ \ \ \ j_2}^{k_{\sigma(1)},k_{\sigma(2)}}\,
\sum_{k_3=1}^n\,G_{k_3,j_3}\,M^{k_3}
-
\delta_{j_1,\ \ \ \ j_3}^{k_{\sigma(1)},k_{\sigma(2)}}\,
\sum_{k_3=1}^n\,G_{k_3,j_2}\,M^{k_3}
+ 
\right. \\
& \
\ \ \ \ \ \ \ \ \ \ \ \ \ \ \
+
\left.
\frac{1}{2}\,\delta_{j_1}^{k_{\sigma(1)}}\,\sum_{k_3=1}^n\, 
H_{k_3,j_3}^{k_{\sigma(2)}}\,L_{j_2}^{k_3} 
-
\frac{1}{2}\,\delta_{j_1}^{k_{\sigma(1)}}\,\sum_{k_3=1}^n\,
H_{k_3,j_2}^{k_{\sigma(2)}}\,L_{j_3}^{k_3}
+ 
\right. \\
& \
\ \ \ \ \ \ \ \ \ \ \ \ \ \ \
\left.
+
\frac{1}{2}\,\delta_{j_2}^{k_{\sigma(1)}}\,\sum_{k_3=1}^n\,
H_{k_3,j_3}^{k_{\sigma(2)}}\,L_{j_1}^{k_3}
-
\frac{1}{2}\,\delta_{j_3}^{k_{\sigma(1)}}\,\sum_{k_3=1}^n\, 
H_{k_3,j_2}^{k_{\sigma(2)}}\,L_{j_1}^{k_3}
+ 
\right. \\
& \
\ \ \ \ \ \ \ \ \ \ \ \ \ \ \
\left.
+
\frac{1}{2}\,\delta_{j_3}^{k_{\sigma(1)}}\,\sum_{k_3=1}^n\,
H_{j_1,j_2}^{k_3}\,L_{k_3}^{k_{\sigma(2)}}
-
\frac{1}{2}\,\delta_{j_2}^{k_{\sigma(1)}}\,\sum_{k_3=1}^n\, 
H_{j_1,j_3}^{k_3}\ L_{k_3}^{k_{\sigma(2)}}
\right).
\endaligned\right.
\end{equation}
\def\theequation{IV'}\begin{equation}
\left\{
\aligned
0
&
=
\sum_{\sigma\in\mathfrak{S}_3}\left( 
\frac{1}{2}\,\delta_{j_3,\ \ \ \ j_1}^{k_{\sigma(3)},k_{\sigma(2)}}\,
L_{j_2,y}^{k_{\sigma(1)}}
- 
\frac{1}{2}\,\delta_{j_2,\ \ \ \ j_1}^{k_{\sigma(3)},k_{\sigma(2)}}\,
L_{j_3,y}^{k_{\sigma(1)}}
+ 
\right. \\
& \ 
\ \ \ \ \ \ \ \ \ \ \ \
\left.
+
\delta_{j_2,\ \ \ \ j_1}^{k_{\sigma(3)},k_{\sigma(2)}}\,
M_{x^{j_3}}^{k_{\sigma(1)}}
- 
\delta_{j_3,\ \ \ \ j_1}^{k_{\sigma(3)},k_{\sigma(2)}}\,
M_{x^{j_2}}^{k_{\sigma(1)}}
+ 
\right. \\
& \
\ \ \ \ \ \ \ \ \ \ \ \
\left.
+
\delta_{j_2,\ \ \ \ j_1}^{k_{\sigma(3)},k_{\sigma(1)}}\, 
\sum_{k_4=1}^n\,H_{k_4,j_3}^{k_{\sigma(2)}}\,M^{k_4}
-
\delta_{j_3,\ \ \ \ j_1}^{k_{\sigma(3)},k_{\sigma(1)}}\, 
\sum_{k_4=1}^n\,H_{k_4,j_2}^{k_{\sigma(2)}}\,M^{k_4}
+ 
\right. \\
& \
\ \ \ \ \ \ \ \ \ \ \ \
\left.
+
\frac{1}{4}\,\delta_{j_1,\ \ \ \ j_3}^{k_{\sigma(1)},k_{\sigma(3)}}\,
\sum_{k_4=1}^n\,L_{k_4}^{k_{\sigma(2)}}\,L_{j_2}^{k_4}
-
\frac{1}{4}\,\delta_{j_1,\ \ \ \ j_2}^{k_{\sigma(1)},k_{\sigma(3)}}\, 
\sum_{k_4=1}^n\,L_{k_4}^{k_{\sigma(2)}}\, L_{j_3}^{k_4}
\right).
\endaligned\right.
\end{equation}
These systems (I'), (II'), (III') and (IV') should be distinguished
from the systems (I), (II), (III) and (IV) of Theorem~1.7 in~\cite{
m2004a}, although they are quite similar. Here, the indices $j_1, j_2,
j_3, k_1, k_2, k_3$ vary in $\{1, 2, \dots, n\}$. By $\mathfrak{
S}_2$ and by $\mathfrak{ S}_3$, we denote the permutation group of
$\{1, 2\}$ and of $\{1, 2, 3\}$. To facilitate hand- and
Latex-writing, partial derivatives are denoted as indices after a
comma; for instance, $G_{j_1, j_2, x^{j_3}}$ is an abreviation for
$\partial G_{j_1, j_2}/ \partial x^{j_3}$. To deduce (I'), (II'),
(III') and (IV') from equation~\thetag{ 
1.25} of~\cite{ m2004a}, it suffices to
observe that every cubic polynomial equation of the form
\def\theequation{1.6}\begin{equation}
\aligned
0 
&
\equiv
A 
+
\sum_{k_1=1}^n\,B_{k_1}\cdot\,
y_{x^{k_1}}
+
\sum_{k_1=1}^n\,\sum_{k_2=1}^n\,C_{k_1,k_2}\cdot\, 
y_{x^{k_1}}\,y_{x^{k_2}}
+ \\
& \ 
\ \ \ \ \ \ \ 
+
\sum_{k_1=1}^n\,\sum_{k_2=1}^n\,\sum_{k_2=1}^n\,D_{k_1,k_2,k_3}\cdot\,
y_{x^{k_1}}\,y_{x^{k_2}}\,y_{x^{k_3}}
\endaligned
\end{equation}
(as for instance~\thetag{ 1.25} in~\cite{ m2004a}) is equivalent to the
annihilation of the following symmetric sums of its coefficients:
\def\theequation{1.7}\begin{equation}
\left\{
\aligned
0 
&
=
A, 
\\
0 
&
=
B_{k_1}, 
\\
0 
&
=
C_{k_1,k_2}
+
C_{k_2,k_1}, 
\\ 
0 
&
=
D_{k_1,k_2,k_3}
+
D_{k_3,k_1,k_2}
+
D_{k_2,k_3,k_1}
+
D_{k_2,k_1,k_3}
+
D_{k_3,k_2,k_1}
+
D_{k_1,k_3,k_2}.
\endaligned\right.
\end{equation}
for all $k_1, k_2, k_3 = 1, \dots, n$.

In conclusion, the functions $G_{j_1, j_2}$, $H_{j_1, j_2}^{k_1}$,
$L_{j_1}^{k_1}$ and $M^{k_1}$ in the statement of Theorem~1.1 are far
from being arbitrary: they satisfy the complicated system of first
order partial differential equations (I'), (II'), (III') and (IV')
above. 

Our proof of Theorem~1.1 is similar to the one provided in~\cite{
m2004a}, in the case of systems of Newtonian particles, so that
most steps of the proof will be summarized. In the end of
this paper, we will delineate a complicated system of {\it second}\,
order partial differential equations satisfied by $G_{j_1, j_2}$,
$H_{j_1, j_2}^{k_1}$, $L_{j_1}^{k_1}$ and $M^{k_1}$ which is the exact
analog of the system described in the abstract. The main technical
part of the proof of Theorem~1.1 will be to establish that this second
order system is a consequence, by linear combinations and by
differentiations, of the first order system (I'), (II'), (III') and
(IV'). Before proceeding further, let us describe rapidly a second 
proof of Theorem~1.1, confirming its validity.

\subsection*{1.8.~Confirmation of Theorem~1.1 
by the method of equivalence} 
Let us summarize the strategy of~\cite{ bn2002}, \cite{ n2003}, 
\cite{ bi2003}. In~\cite{ch1975}, inspired by M.~Hachtroudi's
thesis~\cite{ ha1937}, S.-S.~Chern conducted the equivalence algorithm
associated to the system~\thetag{ 1.2}. His computations are
essentially equivalent to M.~Hachtroudi's (though with much less
information), but they are done in a {\it non-parametric}\, way. It
was deeply known to \'Elie Cartan and it is now well known by modern
followers ({\it cf.} for instance~\cite{ g1989}, \cite{ gtw1989},
\cite{ hk1989}, \cite{ fe1995}, \cite{ ol1995}) that achieving
parametric computations in an equivalence problem is an extremely hard
task, often devoted to computers nowadays. In fact, the so-called {\sl
Cartan Lemma} (\cite{ ol1995}, p.~26) was essentially devised by
\'E.~Cartan as a shortcut, alongside his groundbreaking progresses,
after some experiences of explicit hand computations. The trick in
this lemma is to bypass the (often too long) parametric computations,
nevertheless keeping track of some important informations. Nowadays,
to study an equivalence problem, one often achieves the non-parametric
computations by hand, leaving the explorations of parametric
computations to a computer. In this respect, the thesis~\cite{ n2003}
of Sylvain Neut is extremely interesting, since it implements
the general equivalence algorithm as a Maple package.

However, computer programs achieving formal algebrico-differential
computations {\it with a general number $n$ of variables}\, are still
undevelopped. Consequently, it is interesting to achieve the computations
of~\cite{ ch1975} in a parametric way, because at the end, the
vanishing of all the invariant tensors that one obtains in the final
$\{e\}$-structure would show explicitely under which condition the
system~\thetag{ 1.2} is equivalent to the system $Y_{ X^{ j_1} X^{j_2}} =
0$. To our knowledge, this problem is considered as open in the field
of {\sl symmetries of Cauchy-Riemann submanifolds of $\C^n$}, a
subfield of the mathematics area called {\sl Several Complex
Variables}. However, no specialist of the subfield seems to be aware
of the old reference~\cite{ ha1937}, a text only read, apparently,
by~S.-S.~Chern, who was a student of \'E.~Cartan contemporary to
M.~Hachtroudi.

Let us expose briefly how one obtains the final $\{ e\}$-structure
attached to the equivalence problem associated to the system~\thetag{
1.2}.  Since we have not been able to follow everything in the ancient
text~\cite{ ha1937}, we follow~\cite{ ch1975} with the only change
that we do not introduce the imaginary number $i = \sqrt{ -1}$ in our
structure equations and also, we make some minor changes of sign;
in~\cite{ ch1975}, the computations are conducted over $\K = \C$ with
the idea that they apply to Complex Analysis, but they do hold as well
over $\K = \R$, or even over an arbitary commutative field of
characteristic zero equipped with a valuation (in order to provide
$\K$-analytic functions).

To begin with, consider the following family of $2n+1$ initial
differential forms:
\def\theequation{1.9}\begin{equation}
\left\{
\aligned
\widetilde{\omega}
&
:= 
dy 
- 
\sum_\beta\,y_{x^\beta}\, dx^\beta, 
\\
\widetilde{\omega}^\alpha
&
:=
dx^\alpha, 
\\
\widetilde{\omega}_\alpha
&
:=
dy_{x^\alpha}
-
\sum_\beta\,F^{\alpha,\beta}\,dx^\beta.
\endaligned\right.
\end{equation}
Here, the Greek indices $\alpha$, $\beta$, $\gamma$, $\delta$,
$\varepsilon$, $\zeta$, {\it etc.}, run from $1$ to $n$. Sums
$\sum_{\alpha = 1}^n$ \, are abbreviated as $\sum_\alpha$. In the
paragraph preceding the statement of Theorem~1.7 in~\cite{ m2004a}, we
explain why we cannot use coherently the Einstein summation convention
throughout the paper. For the same reason, we shall abandon this
convention in the present paper. However, we would like to mention
that in~\thetag{ 1.9} above and in \thetag{ 1.10}, \thetag{ 1.11},
\thetag{ 1.12}, \thetag{ 1.13}, \thetag{ 1.14} and~\thetag{ 1.15}
below, the summation convention applies coherently, so that the reader
may drop the sums if (s)he is used to.

Following~\cite{ ch1975}, define the local Lie
group consisting of $(2n+1)\times (2n+1)$ matrices of the form
\def\theequation{1.10}\begin{equation}
g
:=
\left(
\begin{array}{ccc}
u & 0 & 0 \\
u^\alpha & u_\beta^\alpha & 0 \\
u_\alpha & 0 & u\,{u'}_\alpha^\beta
\end{array}
\right),
\end{equation} 
where $u^\alpha$ and $u_\alpha$ are $n \times 1$ vectors close to the
zero vector, where $u_\beta^\alpha$ is a $n \times n$ invertible
matrix close to the identity matrix, where $u$ is a scalar close to
$1$ and where ${u' }_\beta^\alpha$ denotes the inverse matrix of
$u_\beta^\alpha$. With this group, define the {\sl initial
$G$-structure}, consisting of the following collection of $(2n +1)$
differential forms, called the {\sl lifted coframe}:
\def\theequation{1.11}\begin{equation}
\left\{
\aligned
\omega
&
:=
u\cdot\widetilde{\omega},
\\
\omega^\alpha
&
:=
u^\alpha\cdot\widetilde{\omega}
+
\sum_\beta\,u_\beta^\alpha
\cdot
\widetilde{\omega}^\beta,
\\
\omega_\alpha
&
:=
u_\alpha\cdot\widetilde{\omega}
+
\sum_\beta\,u\,{u'}_\alpha^\beta
\cdot
\widetilde{\omega}_\beta.
\endaligned\right.
\end{equation} 
These forms depend both on the ``horizontal'' variables $(x^i, y)$ and
on the ``vertical'' (``group'', ``fiber'') variables $u$, $u^\alpha$,
$u_\alpha$ and $u_\beta^\alpha$. The introduction of this lifted
coframe may be justified as follows.

It is the very first step
of the method of equivalence 
to check that there exists a transformation $(x^i, y)
\mapsto ( \overline{ x}^i, \overline{ y })$ of the completely
integrable system $y_{x^\alpha x^\beta} = F^{\alpha, \beta}$ in the
coordinates $(x^i, y)$ to another completely integrable system
$\overline{ y }_{ \overline{ x}^\alpha \overline{ x}^\beta} =
\overline{ F}^{ \alpha, \beta}$ in other coordinates $( \overline{
x}^i, \overline{ y})$ {\it if and only if}\, there exist functions
$\widehat{ u}$, $\widehat{ u}^\alpha$, $\widehat{ u}_\alpha$,
$\widehat{ u}_\beta^\alpha$ and $\widehat{ v}_\beta^\alpha$ of $(x^i,
y)$, which depend on the functions $(
\overline{ x}^i, \overline{ y })$, on their partial
derivatives with respect to $(x^i, y)$,
on the $y_{x^l}$ and on the $F^{\gamma, \delta}$, such
that the following matrix identity holds:
\def\theequation{1.12}\begin{equation}
\left(
\begin{array}{c}
\overline{\widetilde{\omega}} \\ 
\overline{\widetilde{\omega}}^\alpha \\ 
\overline{\widetilde{\omega }}_\alpha
\end{array}
\right)
=
\left(
\begin{array}{ccc}
\widehat{u} & 0 & 0 \\
\widehat{u}^\alpha & \widehat{u}_\beta^\alpha & 0 \\
\widehat{u}_\alpha & 0 & \widehat{v}_\alpha^\beta
\end{array}
\right)
\left(
\begin{array}{c}
\widetilde{\omega} \\ 
\widetilde{\omega}^\beta \\ 
\widetilde{\omega}_\beta
\end{array}
\right),
\end{equation} 
where the $(2n +1)$ differential forms $\left( \overline{ \widetilde{
\omega }}, \overline{ \widetilde{ \omega }}^\alpha, \overline{
\widetilde{ \omega }}_\alpha \right)$ are defined similarly as
in~\thetag{ 1.9}, in the barred coordinates. Of course, it is
understood that in the left hand side of this matrix identity, the
variables $(\overline{ x}^i, \overline{ y})$ are replaced by their
values with respect to the variables $(x^i, y)$ (in the language of
modern differential geometry, one usually speaks of ``pull-back''). By
a more careful examination of the explicit expressions of the
functions $\widehat{ u }$, $\widehat{ u }^\alpha$, $\widehat{ u
}_\alpha$, $\widehat{ u }_\beta^\alpha$ and $\widehat{ v
}_\beta^\alpha$ in terms of the $F^{\gamma, \delta}$, one observes
that $\widehat{ v}_\alpha^\beta = \widehat{ u} \, \widehat{
u'}_\alpha^\beta$ (\cite{ m2003}). Thus, {\it the collection of
differential forms $(\omega, \omega^\alpha, \omega_\alpha)$ is defined
modulo multiplication by a matrix of functions of $(x^i, y)$ which is
of the specific form~\thetag{ 1.10}}. Based on this preliminary
observation, the general procedure of the equivalence method
(\cite{ g1989}, \cite{ ol1995})
associates to the system~\thetag{ 1.2} 
the lifted coframe~\thetag{ 1.11}.

Applying the exterior differential operator $d$ to $\omega$, to
$\omega^\alpha$ and to $\omega_\alpha$ and absorbing the torsion, it
is shown in~\cite{ ch1975} that the initial structure equations
may be written under the form
\def\theequation{1.13}\begin{equation}
\left\{
\aligned
d\omega
&
=
\sum_\alpha\,\omega^\alpha\wedge\omega_\alpha
+
\omega\wedge\varphi,
\\
d\omega^\alpha
&
=
\sum_\beta\,\omega^\beta\wedge\varphi_\beta^\alpha
+
\omega\wedge\varphi^\alpha, 
\\
d\omega_\alpha
&
=
\sum_\beta\,\varphi_\alpha^\beta\wedge\omega_\beta
+
\omega_\alpha\wedge\varphi
+
\omega\wedge\varphi_\alpha,
\endaligned\right.
\end{equation}
where $\varphi$, $\varphi_\beta^\alpha$, $\varphi^\alpha$ and
$\varphi_\alpha$ are modified Maurer-Cartan forms. We notice that
there are no torsion coefficient to normalize. In fact, the choice of
$v_\beta^\alpha := u \, {u'}_\beta^\alpha$ made in advance above
corresponds to having achieved a first normalization implicitely.

Since the dimension of the Lie symmetry group of the system~\thetag{
1.2} is always finite and in fact bounded by $n^2+ 4n +3$ (\cite{
ha1937}, \cite{ cm1974}, \cite{ ch1975}, \cite{ su2001}, \cite{
gm2003}), according to the general procedure of the method of
equivalence (\cite{ g1989}, \cite{ ol1995}), one simply has to prolong
the initial $G$-structure. One could also apply the \'Elie Cartan
involutivity test to deduce that it is necessary to prolong.

Now, we summarize the remainder of~\cite{ ch1975} very rapidly. After
one prolongation, two normalizations and one supplementary
prolongation, the final $\{e\}$-structure is of the following form:
\def\theequation{1.14}\begin{equation}
\small
\left\{
\aligned
d\omega
&
=
\sum_\alpha\,\omega^\alpha\wedge\omega_\alpha
+
\omega\wedge\varphi,
\\
d\omega^\alpha
&
=
\sum_\beta\,\omega^\beta\wedge\varphi_\beta^\alpha
+
\omega\wedge\varphi^\alpha,
\\
d\omega_\alpha
&
=
\sum_\beta\,\varphi_\alpha^\beta\wedge\omega_\beta
+
\omega_\alpha\wedge\varphi
+
\omega\wedge\varphi_\alpha, 
\\
d\varphi
&
=
\sum_\alpha\,\omega^\alpha\wedge\varphi_\alpha
-
\sum_\alpha\,\omega_\alpha\wedge\varphi^\alpha
+
\omega\wedge\psi,
\\
d\varphi_\beta^\alpha
&
=
\sum_\rho\,\sum_\sigma\,S_{\beta\rho}^{\alpha\sigma}
\cdot
\omega^\rho\wedge\omega_\sigma
+
\sum_\gamma\,R_{\beta\gamma}^\alpha
\cdot
\omega\wedge\omega^\gamma
+
\sum_\gamma\,T_\beta^{\alpha\gamma}
\cdot
\omega\wedge\omega_\gamma
+ \\
& \
\ \ \ \ \
+
\omega^\alpha\wedge\varphi_\beta
+
\delta_\beta^\alpha\,\sum_\gamma\,\omega^\gamma\wedge\varphi_\gamma
+
\omega_\beta\wedge\varphi^\alpha
+
\sum_\gamma\,\varphi_\beta^\gamma\wedge\varphi_\gamma^\alpha
+
\frac{1}{2}\,\delta_\beta^\alpha
\cdot
\omega\wedge\psi, 
\\
d\varphi^\alpha
&
=
\sum_\beta\,\sum_\gamma\,T_\beta^{\alpha\gamma}
\cdot
\omega^\beta\wedge\omega_\gamma
+
\frac{1}{2}\,\sum_\beta\,Q_\beta^\alpha
\cdot
\omega\wedge\omega^\beta
+
\sum_\beta\,L^{\alpha\beta}
\cdot
\omega\wedge\omega_\beta
+ \\
& \
\ \ \ \ \
+
\varphi\wedge\varphi^\alpha
-
\sum_\beta\,\varphi_\beta^\alpha\wedge\varphi^\beta
+
\frac{1}{2}\,\omega^\alpha\wedge\psi, 
\\
d\varphi_\alpha
&
=
\sum_\beta\,\sum_\gamma\,R_{\alpha\beta}^\gamma
\cdot
\omega^\beta\wedge\omega_\gamma
+
\frac{1}{2}\,\sum_\beta\,Q_\alpha^\beta
\cdot
\omega\wedge\omega_\beta
+
\sum_\beta\,P_{\alpha\beta}
\cdot
\omega\wedge\omega^\beta
+ \\
& \
\ \ \ \ \
+
\sum_\beta\,\varphi_\alpha^\beta\wedge\varphi_\beta
+
\frac{1}{2}\,\omega_\alpha\wedge\psi, 
\\
d\psi
&
=
\sum_\alpha\,\sum_\beta\,Q_\alpha^\beta
\cdot\omega^\alpha\wedge\omega_\beta
+
\sum_\alpha\,H_\alpha
\cdot
\omega\wedge\omega^\alpha
+
\sum_\alpha\,K^\alpha
\cdot
\omega\wedge\omega_\alpha
+ \\
& \
\ \ \ \ \ 
+
\varphi\wedge\psi
+
2\,\sum_\alpha\,\varphi^\alpha\wedge\varphi_\alpha.
\endaligned\right.
\end{equation}
These structure equations incorporate $8$ families $\omega$,
$\omega^\alpha$, $\omega_\alpha$, $\varphi$, $\varphi_\beta^\alpha$,
$\varphi^\alpha$, $\varphi_\alpha$, $\psi$ of differential forms, of
total cardinality $n^2 + 4n + 3$, together with $8$ invariant tensors
$S_{\beta \rho}^{\alpha \sigma}$, $R_{\beta \gamma}^\alpha$,
$T_\beta^{\alpha \gamma}$, $Q_\beta^\alpha$, $L^{\alpha \beta}$,
$P_{\alpha \beta}$, $H_\alpha$ and $K_\alpha$ having some specific
index symmetries that we shall not use.

Applying the exterior differential operator $d$ to these $8$ families
of structure equations~\thetag{ 1.14}, one verifies that the seven
invariant tensors $R_{\beta \gamma}^\alpha$, $T_\beta^{\alpha
\gamma}$, $Q_\beta^\alpha$, $L^{\alpha \beta}$, $P_{\alpha \beta}$,
$H_\alpha$ and $K_\alpha$ are in fact functionally dependent on the
{\sl fundamental tensors}\, $S_{\beta \rho}^{\alpha \sigma}$, namely
they are certain coframe derivatives of the $S_{\beta \rho}^{\alpha
\sigma}$. In the (very similar) context of the equivalence problem
associated with a Levi non-degenerate hypersurface of $\C^{n+1}$, this
computation was achieved by S.M.~Webster in the Appendix of~\cite{
cm1974}; in the precise context of the equivalence problem associated
to the system~\thetag{ 1.2}, this computation is not achieved
in~\cite{ ch1975}, but {\it see} \cite{ bn2002}, \cite{ bi2003},
\cite{ n2003} and~\cite{ m2003} for details. 
For us, the precise nature of this
functional dependence does not matter.

In fact, it is only in the computer science thesis~\cite{ n2003} that
the complete explicit parametric computation of the $8$ families
of differential forms
$\omega$, $\omega^\alpha$, $\omega_\alpha$, $\varphi$,
$\varphi_\beta^\alpha$, $\varphi^\alpha$, $\varphi_\alpha$, $\psi$
together with the $8$ invariant tensors $S_{\beta \rho}^{\alpha
\sigma}$, $R_{\beta \gamma}^\alpha$, $T_\beta^{\alpha \gamma}$,
$Q_\beta^\alpha$, $L^{\alpha \beta}$, $P_{\alpha \beta}$, $H_\alpha$
and $K_\alpha$ is achieved, in the case $n=2$ and with the help of
Maple. Although the task is really of impressive size, the author of
this paper has the project of achieving manually the general
computation for $n\geq 2$ variables (\cite{ m2005}). In~\cite{
ha1937}, the parametric computations are achieved completely only in
the case $n=2$. According to the author's experience, it appears that
the tensorial formalism helps to shorten importantly the size of the
electronic computations and that the case $n\geq 2$ is not much more
difficult than the case $n=2$, as one may already observe by
reading~\cite{ m2004a}. This is why the project~\cite{ m2005} seems to
be an accessible task. A similar project would be to compute in a
parametric way the structure equations obtained in~\cite{ fe1995} for
the equivalence problem associated to a system of Newtonian particles;
the computations are quite similar, indeed, as well as the
computations of this paper as quite similar to the computations
of~\cite{ m2004a}. This will be achieved in the future, if
time permits.

At present, fortunately, there is a kind of ``miracle'': the
interesting tensor $S_{\beta \rho}^{\alpha \sigma}$ appears
essentially at the {\it beginning}\, of the computation of the final
$\{e\}$-structure. Thus, it is not necessary to go up to the end of
the algorithm in order to deduce the final parametric expression of $S_{\beta
\rho}^{\alpha \sigma}$, and especially to characterize the maximally
symmetric systems~\thetag{ 1.2}, {\it i.e.} those for which all
the $8$ tensors $S_{\beta \rho}^{\alpha \sigma}$, $R_{\beta
\gamma}^\alpha$, $T_\beta^{\alpha \gamma}$, $Q_\beta^\alpha$,
$L^{\alpha \beta}$, $P_{\alpha \beta}$, $H_\alpha$ and $K_\alpha$
vanish. In~\cite{ m2003}, we obtained:
\def\theequation{1.15}\begin{equation}
\left\{
\aligned
S_{\beta\rho}^{\alpha\sigma}
&
=
\delta_\rho^\sigma
\left(
\frac{1}{n+2}\,\sum_\gamma\,\sum_\delta\,\sum_\varepsilon\,
u^{-1}\,{u'}_\beta^\delta\,u_\varepsilon^\alpha\,
F_{y_{x^\gamma}y_{x^\varepsilon}}^{\gamma,\delta}
\right)
+ \\
& \
\ \ \ \ \
+
\delta_\rho^\alpha
\left(
\frac{1}{n+2}\,\sum_\gamma\,\sum_\delta\,\sum_\varepsilon\,
u^{-1}\,{u'}_\beta^\delta\,u_\varepsilon^\sigma\,
F_{y_{x^\gamma}y_{x^\varepsilon}}^{\gamma,\delta}
\right)
+ \\
& \
\ \ \ \ \
+
\delta_\beta^\sigma
\left(
\frac{1}{n+2}\,\sum_\gamma\,\sum_\delta\,\sum_\varepsilon\,
u^{-1}\,{u'}_\rho^\delta\,u_\varepsilon^\alpha\,
F_{y_{x^\gamma}y_{x^\varepsilon}}^{\gamma,\delta}
\right)
+ \\
& \
\ \ \ \ \
+
\delta_\beta^\alpha
\left(
\frac{1}{n+2}\,\sum_\gamma\,\sum_\delta\,\sum_\varepsilon\,
u^{-1}\,{u'}_\rho^\delta\,u_\varepsilon^\sigma\,
F_{y_{x^\gamma}y_{x^\varepsilon}}^{\gamma,\delta}
\right)
- \\
& \
\ \ \ \ \
-
\left(
\delta_\rho^\sigma\,\delta_\beta^\alpha
+
\delta_\rho^\alpha\,\delta_\beta^\sigma
\right)
\cdot
\left(
\frac{1}{(n+1)(n+2)}\,\sum_\gamma\,\sum_\delta\,
u^{-1}\,F_{y_{x^\gamma}y_{x^\delta}}^{\gamma,\delta}
\right)
- \\
& \
\ \ \ \ \
-
\sum_\gamma\,\sum_\delta\,\sum_\varepsilon\,\sum_\zeta\,
u^{-1}\,{u'}_\rho^\delta\,u_\zeta^\sigma\,{u'}_\beta^\gamma\,
u_\varepsilon^\alpha\,F_{y_{x^\varepsilon}y_{x^\zeta}}^{\gamma,\delta}.
\endaligned\right.
\end{equation}
It may be verified ({\it cf.} also~\cite{ bi2003}) that the general
tensor $S_{\beta \rho }^{ \alpha \sigma}$ is a certain ``tensorial
rotation'' of its value ``at the identity'', namely its value when the
matrix~\thetag{ 1.10} is the identity; it suffices to put $u:= 1$,
$u_\beta^\alpha := \delta_\beta^\alpha$ and ${u'}_\beta^\alpha :=
\delta_\beta^\alpha$ in~\thetag{ 1.15} above. Next, it may be
verified that the vanishing of the value of $S_{\beta \rho }^{ \alpha
\sigma}$ at the identity (which is equivalent to the vanishing of the
general $S_{\beta \rho }^{ \alpha \sigma}$) provides a system of
second linear second order partial differential equations involving
only the partial derivatives $F_{y_{ x^\varepsilon }y_{ x^\zeta }}^{
\gamma, \delta}$. Finally, one verifies (\cite{ bi2003}) that this
last system is equivalent to the fact that the $F^{ j_1, j_2}$ are of
the specific cubic polynomial form~\thetag{ 1.3}. In conclusion, this
provides a second (very summarized) proof of Theorem~1.1.

In the remainder of the paper, we shall not speak anymore of the
equivalence method and we will focus on describing precisely how
to establish Theorem~1.1 using S.~Lie's original techniques, as
in~\cite{ m2004a}.

\subsection*{1.17.~Acknowledgment}
We are very grateful to Sylvain Neut and to Michel Petitot (LIFL,
University of Lille~1), who implemented the \'Elie Cartan equivalence
algorithm and who discovered the validity of Theorem~1.1 in the case
$n=2$, and also of Theorem~1.7 of~\cite{ m2004a} in the case
$m=2$. Without these discoveries, we would not have pushed our manual
computations up to the very end of the proof of Theorem~1.1 (and of
Theorem~1.7 in~\cite{ m2004a}). We also acknowledge interesting
exchanges about the method of equivalence
with Camille Bi\`eche, Sylvain Neut and Michel Petitot.

\subsection*{1.18.~Closing remark}
After Theorem~1.1 was established, namely after the references~\cite{
bn2002}, \cite{ n2003}, \cite{ bi2003} and \cite{ m2003} were
completed, we discovered at the mathematical library of the \'Ecole
Normale Sup\'erieure of Paris that Mohsen Hachtroudi, an Iranian
student of \'Elie Cartan, also obtained a proof of Theorem~1.1 (\cite{
ha1937}, p.~53), sixty seven years ago. However, his computations
follow the method of equivalence, in the spirit of his master, whereas
we conduct ours in the spirit of Sophus Lie. Also, in~\cite{ ha1937},
M.~Hachtroudi only refers to the abbreviated formulation~\thetag{ 1.4}
of the compatibility conditions, and one does not find there the
explicit and complete formulation of the systems (I'), (II'), (III')
and (IV'). Finally, the combinatorial formulas that we provide
in~\thetag{ 2.10}, in~\thetag{ 2.14} in~\thetag{ 2.15}, in~\thetag{
2.22}, in~\thetag{ 2.23}, 
in~\thetag{ 2.31} and in~\thetag{ 2.35} below seem to be new in
the Lie theory of symmetries of partial differential equations ({\it
see} also~\cite{ m2004b}).

\section*{ \S2.~Completely integrable systems of \\
second order ordinary differential equations}

\subsection*{ 2.1.~Prolongation of a point transformation to 
the second order jet space} Let $\K = \R$ or $\C$, let
$n\in \N$, 
\underline{suppose
$n\geq 2$}, let $x = (x^1, \dots, x^n) \in \K^n$ and let $y
\in \K$. According to the main assumption of Theorem~1.1, we have to
consider a local $\K$-analytic diffeomorphism of the form
\def\theequation{2.2}\begin{equation}
\left(
x^{j_1}, y 
\right)
\longmapsto 
\left(
 X^j (x^{j_1}, y), 
\
Y( x^{j_1}, y)
\right),
\end{equation}
which transforms the system~\thetag{ 1.2} to the system $Y_{X^{i_1}
X^{i_2}}= 0$, $1\leq j_1, j_2 \leq n$. Without loss of generality, we
shall assume that this transformation is close to the identity. To
obtain the precise expression~\thetag{ 2.35} 
of the transformed system~\thetag{ 1.2}, we have to prolong the above
diffeomorphism to the second order jet space. We introduce the
coordinates $\left( x^j, y, y_{x^{j_1}}, y_{x^{j_1}x^{j_2}} \right)$
on the second order jet space. Let
\def\theequation{2.3}\begin{equation}
D_k
:= 
\frac{\partial}{\partial x^k}
+
y_{x^k}\,\frac{\partial}{\partial y}
+
\sum_{l=1}^n\,y_{x^kx^l}\,\frac{\partial}{\partial y_{x^l}},
\end{equation}
be the $k$-th total differentiation operator. According to~\cite{
bk1989}, for the first order partial derivatives, one has the
(implicit, compact) expression:
\def\theequation{2.4}\begin{equation}
\left(
\begin{array}{c}
Y_{X^1} \\
\vdots \\
Y_{X^n}
\end{array}
\right)
=
\left(
\begin{array}{ccc}
D_1X^1 & \cdots & D_1X^n \\
\vdots & \cdots & \vdots \\
D_nX^1 & \cdots & D_nX^n
\end{array}
\right)^{-1} 
\left(
\begin{array}{c}
D_1Y \\
\vdots \\
D_nY
\end{array}
\right),
\end{equation}
where $(\cdot)^{ -1}$ denotes the inverse matrix, which exists, since
the transformation~\thetag{ 2.2} is close to the identity. For the
second order partial derivatives, again according to~\cite{ bk1989},
one has the (implicit, compact) expressions:
\def\theequation{2.5}\begin{equation}
\left(
\begin{array}{c}
Y_{X^jX^1} \\
\vdots \\
Y_{X^jX^n}
\end{array}
\right)
=
\left(
\begin{array}{ccc}
D_1X^1 & \cdots & D_1X^n \\
\vdots & \cdots & \vdots \\
D_nX^1 & \cdots & D_nX^n
\end{array}
\right)^{-1} 
\left(
\begin{array}{c}
D_1Y_{X^j} \\
\vdots \\
D_nY_{X^j}
\end{array}
\right),
\end{equation}
for $j = 1, \dots, n$. Let $DX$ denote the matrix $\left( D_iX^j
\right)_{ 1\leq i\leq n}^{ 1\leq j\leq n}$, where $i$ is the index of
lines and $j$ the index of columns, let $Y_X$ denote the column matrix
$\left( Y_{X^i} \right)_{ 1\leq i \leq n}$ and let $DY$ be the column
matrix $\left( D_i Y \right)_{ 1\leq i \leq n}$.

By inspecting~\thetag{ 2.5} above, we see that the
equivalence between {\bf (i)}, {\bf (ii)}
and {\bf (iii)} just below is obvious:

\def\thelemma{2.6}\begin{lemma}
The following conditions are equivalent{\rm :}
\begin{itemize}
\item[{\bf (i)}]
the differential equations $Y_{X^jX^k} = 0$ hold for $1\leq j,k \leq
n${\rm ;}
\item[{\bf (ii)}]
the matrix equations $D_k( Y_X)= 0$ hold for $1\leq k\leq n${\rm ;}
\item[{\bf (iii)}]
the matrix equations $DX \cdot D_k( Y_X)= 0$ hold for $1\leq 
k\leq n${\rm ;}
\item[{\bf (iv)}]
the matrix equations
$0 = D_k( DX) \cdot Y_X - D_k (DY)$ hold for $1\leq k\leq n$.
\end{itemize}
\end{lemma}

Formally, in the sequel, it will be more convenient to achieve the
explicit computations starting from condition {\bf (iv)}, since no
matrix inversion at all is involved in it.

\proof
Indeed, applying the
total differentiation operator $D_k$ to 
the matrix equation~\thetag{ 2.4} written under the equivalent
form $0 = DX\cdot Y_X - DY$, we get:
\def\theequation{2.7}\begin{equation}
0 
=
D_k(DX)\cdot Y_X 
+
DX\cdot D_x(Y_X)
-
D_k(DY),
\end{equation}
so that the equivalence between {\bf (iii)} and
{\bf (iv)} is now clear.
\endproof

\subsection*{ 2.8.~An explicit formula in the case $n = 2$}
Thus, we can start to develope explicitely the 
matrix equations
\def\theequation{2.9}\begin{equation}
0 
=
D_k (DX) \cdot Y_X 
- 
D_k(DY).
\end{equation}
In it, some huge formal expressions are hidden behind the symbol
$D_k$. Proceeding inductively, we start by examinating the case $n=2$
thoroughly. By direct computations which require to be clever, we
reconstitute some $3\times 3$ determinants in the four (in 
fact three) developed
equations~\thetag{ 2.9}. After some work, the first equation is:
\def\theequation{2.10}\begin{equation}
\aligned
0
& 
= 
y_{x^1x^1}
\cdot
\left\vert
\begin{array}{ccc}
X_{x^1}^1 & X_{x^2}^1 & X_y^1 \\
X_{x^1}^2 & X_{x^2}^2 & X_y^2 \\
Y_{x^1} & Y_{x^2} & Y_y \\
\end{array}
\right\vert
+
\left\vert
\begin{array}{ccc}
X_{x^1}^1 & X_{x^2}^1 & X_{x^1x^1}^1 \\
X_{x^1}^2 & X_{x^2}^2 & X_{x^1x^1}^2 \\
Y_{x^1} & Y_{x^2} & Y_{x^1x^1} \\
\end{array}
\right\vert
+ \\
& \
\ \ \ \ \
+
y_{x^1}
\cdot
\left\{
2\,
\left\vert
\begin{array}{ccc}
X_{x^1}^1 & X_{x^2}^1 & X_{x^1y}^1 \\
X_{x^1}^2 & X_{x^2}^2 & X_{x^1y}^2 \\
Y_{x^1} & Y_{x^2} & Y_{x^1y} \\
\end{array}
\right\vert
-
\left\vert
\begin{array}{ccc}
X_{x^1x^1}^1 & X_{x^2}^1 & X_y^1 \\
X_{x^1x^1}^2 & X_{x^2}^2 & X_y^2 \\
Y_{x^1x^1} & Y_{x^2} & Y_y \\
\end{array}
\right\vert
\right\}
+ \\
& \
\ \ \ \ \
+
y_{x^2}
\cdot
\left\{
-
\left\vert
\begin{array}{ccc}
X_{x^1}^1 & X_{x^1x^1}^1 & X_y^1 \\
X_{x^1}^2 & X_{x^1x^1}^2 & X_y^2 \\
Y_{x^1} & Y_{x^1x^1} & Y_y \\
\end{array}
\right\vert
\right\}
+ \\
\endaligned
\end{equation}
$$
\aligned
& \
\ \ \ \ \
+
y_{x^1}\,y_{x^1}
\cdot
\left\{
\left\vert
\begin{array}{ccc}
X_{x^1}^1 & X_{x^2}^1 & X_{yy}^1 \\
X_{x^1}^2 & X_{x^2}^2 & X_{yy}^2 \\
Y_{x^1} & Y_{x^2} & Y_{yy} \\
\end{array}
\right\vert
-
2\,
\left\vert
\begin{array}{ccc}
X_{x^1y}^1 & X_{x^2}^1 & X_y^1 \\
X_{x^1y}^2 & X_{x^2}^2 & X_y^2 \\
Y_{x^1y} & Y_{x^2} & Y_y \\
\end{array}
\right\vert
\right\}
+ \\
& \
\ \ \ \ \ 
+
y_{x^1}\,y_{x^2}
\cdot
\left\{
-2\,
\left\vert
\begin{array}{ccc}
X_{x^1}^1 & X_{x^1y}^1 & X_y^1 \\
X_{x^1}^2 & X_{x^1y}^2 & X_y^2 \\
Y_{x^1} & Y_{x^1y} & Y_y \\
\end{array}
\right\vert
\right\}
+ \\
& \
\ \ \ \ \
+
y_{x^1}\,y_{x^1}\,y_{x^1}
\cdot
\left\{
-
\left\vert
\begin{array}{ccc}
X_{yy}^1 & X_{x^2}^1 & X_y^1 \\
X_{yy}^2 & X_{x^2}^2 & X_y^2 \\
Y_{yy} & Y_{x^2} & Y_y \\
\end{array}
\right\vert
\right\}
+ \\
& \
\ \ \ \ \
+
y_{x^1}\,y_{x^1}\,y_{x^2}
\cdot
\left\{
-
\left\vert
\begin{array}{ccc}
X_{x^1}^1 & X_{yy}^1 & X_y^1 \\
X_{x^1}^2 & X_{yy}^2 & X_y^2 \\
Y_{x^1} & Y_{yy} & Y_y \\
\end{array}
\right\vert
\right\}.
\endaligned
$$
This formula and the two next~\thetag{ 2.22}, 
\thetag{ 2.23} have been checked by Sylvain Neut
and Michel Petitot with the help of Maple.

\subsection*{2.11.~Comparison with the 
coefficients of the second prolongation of a vector field} At present,
it is useful to make an illuminating digression which will help us to
devise what is the general form of the development of the
equations~\thetag{ 2.9}. Consider an arbitrary vector field of the
form
\def\theequation{2.12}\begin{equation}
\mathcal{ L}
:= 
\sum_{k=1}^n\,\mathcal{X}^k\,\frac{\partial}{\partial x^k}
+
\mathcal{Y}\,\frac{\partial}{\partial y},
\end{equation}
where the coefficients $\mathcal{ X}^k$ and $\mathcal{ Y}$ are
functions of $(x^i, y)$. According to~\cite{ ol1986}, \cite{
bk1989}, there exists a unique prolongation $\mathcal{ L}^{(2)}$ of
this vector field to the second order jet space, of the form
\def\theequation{2.13}\begin{equation}
\mathcal{ L}^{(2)} 
:=
\mathcal{L}
+
\sum_{j_1=1}^n\,
{\bf Y}_{j_1}\,\frac{\partial}{\partial y_{x^{j_1}}}
+
\sum_{j_1=1}^n\,\sum_{j_2=1}^n\,
{\bf Y}_{j_1,j_2}\,\frac{\partial}{\partial y_{x^{j_1}x^{j_2}}},
\end{equation}
where the coefficients ${\bf Y}_{ j_1}$, ${\bf Y}_{ j_1, j_2}$ may be
computed by means of some inductive devices explained in~\cite{ ol1986},
\cite{ bk1989}. In~\cite{ gm2003} ({\it see} also~\cite{ su2001} for
formulas which do not involve the Kronecker symbol),
we obtained the following perfect formulas:
\def\theequation{2.14}\begin{equation}
\left\{
\aligned
{\bf Y}_{j_1,j_2}
&
=
\mathcal{Y}_{x^{j_1}x^{j_2}}
+
\sum_{k_1=1}^n\,y_{x^{k_1}}
\cdot
\left\{
\delta_{j_1}^{k_1}\,\mathcal{Y}_{x^{j_2}y}
+
\delta_{j_2}^{k_1}\,\mathcal{Y}_{x^{j_1}y}
-
\mathcal{X}_{x^{j_1}x^{j_2}}^{k_1}
\right\}
+ \\
& \
\ \ \ \ \
+
\sum_{k_1=1}^n\,\sum_{k_2=1}^n\,y_{x^{k_1}}\,y_{x^{k_2}}
\cdot
\left\{
\delta_{j_1,\,j_2}^{k_1,k_2}\,\mathcal{Y}_{yy}
-
\delta_{j_1}^{k_1}\,\mathcal{X}_{x^{j_2}y}^{k_2}
-
\delta_{j_2}^{k_1}\,\mathcal{X}_{x^{j_1}y}^{k_2}
\right\}
+ \\
& \
\ \ \ \ \
+
\sum_{k_1=1}^n\,\sum_{k_2=1}^n\,\sum_{k_3=1}^n\,
y_{x^{k_1}}\,y_{x^{k_2}}\,y_{x^{k_3}}
\cdot
\left\{
-\delta_{j_1,\,j_2}^{k_1,k_2}\,\mathcal{X}_{yy}^{k_3}
\right\},
\endaligned\right.
\end{equation} 
for $j_1, j_2 = 1,\dots, n$.
The expression of ${\bf Y}_{j_1}$ does not
matter for us here.
Specifying this formula to the the case $n=2$ and taking account of
the symmetry ${\bf Y}_{ 1, 2} = {\bf Y}_{ 2,1}$ we get the following
three second order coefficients:
\def\theequation{2.15}\begin{equation}
\left\{
\aligned
{\bf Y}_{1,1}
&
=
\mathcal{Y}_{x^1x^1}
+
y_{x^1}
\cdot
\left\{
2\,\mathcal{Y}_{x^1y}
-
\mathcal{X}_{x^1x^1}^1
\right\}
+
y_{x^2}
\cdot
\left\{
-
\mathcal{X}_{x^1x^1}^2
\right\}
+ \\
& \
\ \ \ \ \
+
y_{x^1}\,y_{x^1}
\cdot
\left\{
\mathcal{Y}_{yy}
-
2\,\mathcal{X}_{x^1y}^1
\right\}
+
y_{x^1}\,y_{x^2}
\cdot
\left\{
-
2\,\mathcal{X}_{x^1y}^2
\right\}
+ \\
& \
\ \ \ \ \ 
+
y_{x^1}\,y_{x^1}\,y_{x^1}
\cdot
\left\{
-
\mathcal{X}_{yy}^1
\right\}
+
y_{x^1}\,y_{x^1}\,y_{x^2}
\cdot
\left\{
-
\mathcal{X}_{yy}^2
\right\}, \\
{\bf Y}_{1,2}
&
=
\mathcal{Y}_{x^1x^2}
+
y_{x^1}
\cdot
\left\{
\mathcal{Y}_{x^2y}
-
\mathcal{X}_{x^1x^2}^1
\right\}
+
y_{x^2}
\cdot
\left\{
\mathcal{Y}_{x^1y}
-
\mathcal{X}_{x^1x^2}^2
\right\}
+ \\
& \
\ \ \ \ \
+
y_{x^1}\,y_{x^1}
\cdot
\left\{
-
\mathcal{X}_{x^2y}^1
\right\}
+
y_{x^1}\,y_{x^2}
\cdot
\left\{
\mathcal{Y}_{yy}
-
\mathcal{X}_{x^1y}^1
-
\mathcal{X}_{x^2y}^2
\right\}
+ \\
& \
\ \ \ \ \
+
y_{x^2}\,y_{x^2}
\cdot
\left\{
-
\mathcal{X}_{x^1y}^2
\right\}
+ \\
& \
\ \ \ \ \ 
+
y_{x^1}\,y_{x^1}\,y_{x^2}
\cdot
\left\{
-
\mathcal{X}_{yy}^1
\right\}
+
y_{x^1}\,y_{x^2}\,y_{x^2}
\cdot
\left\{
-
\mathcal{X}_{yy}^2
\right\}, \\
{\bf Y}_{2,2}
&
=
\mathcal{Y}_{x^2x^2}
+
y_{x^1}
\cdot
\left\{
-
\mathcal{X}_{x^2x^2}^1
\right\}
+
y_{x^2}
\cdot
\left\{
2\,\mathcal{Y}_{x^2y}
-
\mathcal{X}_{x^2x^2}^2
\right\}
+ \\
& \
\ \ \ \ \
+
y_{x^1}\,y_{x^2}
\cdot
\left\{
-
2\,\mathcal{X}_{x^2y}^1
\right\}
+
y_{x^2}\,y_{x^2}
\cdot
\left\{
\mathcal{Y}_{yy}
-
2\,\mathcal{X}_{x^2y}^2
\right\}
+ \\
& \
\ \ \ \ \ 
+
y_{x^1}\,y_{x^2}\,y_{x^2}
\cdot
\left\{
-
\mathcal{X}_{yy}^1
\right\}
+
y_{x^2}\,y_{x^2}\,y_{x^2}
\cdot
\left\{
-
\mathcal{X}_{yy}^2
\right\}.
\endaligned\right.
\end{equation}
We would like to mention that the computation of ${\bf Y}_{j_1, j_2}$,
$1\leq j_1, j_2 \leq 2$, above is easier than the verification
of~\thetag{ 2.10}. Based on the three formulas~\thetag{ 2.15}, we claim
that we can guess the second and the third equations, 
which would be obtained by developing and by simplifying~\thetag{ 2.9},
namely with
$y_{x^1x^2}$ and with $y_{ x^2 x^2}$ instead of $y_{x^1 x^2}$
in~\thetag{ 2.10}. Our dictionary to translate from the first
formula~\thetag{ 2.15} to~\thetag{ 2.10} may be described as
follows. Begin with the {\sl Jacobian determinant}
\def\theequation{2.16}\begin{equation}
\left\vert
\begin{array}{ccc}
X_{x^1}^1 & X_{x^2}^1 & X_y^1 \\
X_{x^1}^2 & X_{x^2}^2 & X_y^2 \\
Y_{x^1} & Y_{x^2} & Y_y \\
\end{array}
\right\vert
\end{equation}
of the change of coordinates~\thetag{ 2.2}. Since this change of
coordinates is close to the identity, we may consider that the following
Jacobian matrix approximation holds:
\def\theequation{2.17}\begin{equation}
\left(
\begin{array}{ccc}
X_{x^1}^1 & X_{x^2}^1 & X_y^1 \\
X_{x^1}^2 & X_{x^2}^2 & X_y^2 \\
Y_{x^1} & Y_{x^2} & Y_y \\
\end{array}
\right)
\cong
\left(
\begin{array}{ccc}
1 & 0 & 0 \\
0 & 1 & 0 \\
0 & 0 & 1\\
\end{array}
\right).
\end{equation}
The jacobian matrix has three columns. There are six possible second
order derivatives with respect to the variables $(x^1, x^2, y)$,
namely
\def\theequation{2.18}\begin{equation}
(\cdot)_{x^1x^1}, 
\ \ \ \ \
(\cdot)_{x^1x^2}, 
\ \ \ \ \
(\cdot)_{x^2x^2}, 
\ \ \ \ \
(\cdot)_{x^1y}, 
\ \ \ \ \
(\cdot)_{x^2y}, 
\ \ \ \ \
(\cdot)_{yy}.
\end{equation}
In the Jacobian determinant~\thetag{ 2.16}, by replacing any one of
the three columns of first order derivatives with a column of second
order derivatives, we obtain exactly $3\times 6 =18$ possible
determinants. For instance, by replacing the third column by the
second order derivative $(\cdot )_{ x^1y}$ or the first column by the
second order derivative $(\cdot )_{ x^1x^1}$, we get:
\def\theequation{2.19}\begin{equation}
\left\vert
\begin{array}{ccc}
X_{x^1}^1 & X_{x^2}^1 & X_{x^1y}^1 \\
X_{x^1}^2 & X_{x^2}^2 & X_{x^1y}^2 \\
Y_{x^1} & Y_{x^2} & Y_{x^1y} \\
\end{array}
\right\vert 
\ \ \ \ \ \ \ \ \ \ 
{\rm or}
\ \ \ \ \ \ \ \ \ \
\left\vert
\begin{array}{ccc}
X_{x^1x^1}^1 & X_{x^2}^1 & X_y^1 \\
X_{x^1x^1}^2 & X_{x^2}^2 & X_y^2 \\
Y_{x^1x^1} & Y_{x^2} & Y_y \\
\end{array}
\right\vert.
\end{equation}
We recover the two determinants appearing in the second line
of~\thetag{ 2.10}. On the other hand, according to the 
approximation~\thetag{ 2.17}, these two determinants are
essentially equal to
\def\theequation{2.20}\begin{equation}
\left\vert
\begin{array}{ccc}
1 & 0 & X_{x^1y}^1 \\
0 & 1 & X_{x^1y}^2 \\
0 & 0 & Y_{x^1y} \\
\end{array}
\right\vert 
=
Y_{x^1y}
\ \ \ \ \ \ \ \ \ \ 
\text{\rm or to}
\ \ \ \ \ \ \ \ \ \
\left\vert
\begin{array}{ccc}
X_{x^1x^1}^1 & 0 & 0 \\
X_{x^1x^1}^2 & 1 & 0 \\
Y_{x^1x^1} & 0 & 1 \\
\end{array}
\right\vert
=
X_{x^1x^1}^1.
\end{equation}
Consequently, in the second line of~\thetag{ 2.10}, up to a change to
calligraphic letters, we recover the coefficient
\def\theequation{2.21}\begin{equation}
2\,\mathcal{Y}_{x^1y}
-
\mathcal{X}_{x^1x^1}^1
\end{equation}
of $y_{x_1}$ in the expression of ${\bf Y}_{1,1}$ in~\thetag{ 2.15}.
In conclusion, we have discovered how to pass symbolically from the
first equation~\thetag{ 2.15} to the equation~\thetag{ 2.10}
and conversely.

Translating the second equation~\thetag{ 2.15}, we deduce, {\it
without any further computation}, that the second equation which would
be obtained by developing~\thetag{ 2.9} in length, is:
\def\theequation{2.22}\begin{equation}
\aligned
0
& 
= 
y_{x^1x^2}
\cdot
\left\vert
\begin{array}{ccc}
X_{x^1}^1 & X_{x^2}^1 & X_y^1 \\
X_{x^1}^2 & X_{x^2}^2 & X_y^2 \\
Y_{x^1} & Y_{x^2} & Y_y \\
\end{array}
\right\vert
+
\left\vert
\begin{array}{ccc}
X_{x^1}^1 & X_{x^2}^1 & X_{x^1x^2}^1 \\
X_{x^1}^2 & X_{x^2}^2 & X_{x^1x^2}^2 \\
Y_{x^1} & Y_{x^2} & Y_{x^1x^2} \\
\end{array}
\right\vert
+ \\
& \
\ \ \ \ \
+
y_{x^1}
\cdot
\left\{
\left\vert
\begin{array}{ccc}
X_{x^1}^1 & X_{x^2}^1 & X_{x^2y}^1 \\
X_{x^1}^2 & X_{x^2}^2 & X_{x^2y}^2 \\
Y_{x^1} & Y_{x^2} & Y_{x^2y} \\
\end{array}
\right\vert
-
\left\vert
\begin{array}{ccc}
X_{x^1x^2}^1 & X_{x^2}^1 & X_y^1 \\
X_{x^1x^2}^2 & X_{x^2}^2 & X_y^2 \\
Y_{x^1x^2} & Y_{x^2} & Y_y \\
\end{array}
\right\vert
\right\}
+ \\
& \
\ \ \ \ \
+
y_{x^2}
\cdot
\left\{
\left\vert
\begin{array}{ccc}
X_{x^1}^1 & X_{x^2}^1 & X_{x^1y}^1 \\
X_{x^1}^2 & X_{x^2}^2 & X_{x^1y}^2 \\
Y_{x^1} & Y_{x^2} & Y_{x^1y} \\
\end{array}
\right\vert
-
\left\vert
\begin{array}{ccc}
X_{x^1}^1 & X_{x^1x^2}^1 & X_y^1 \\
X_{x^1}^2 & X_{x^1x^2}^2 & X_y^2 \\
Y_{x^1} & Y_{x^1x^2} & Y_y \\
\end{array}
\right\vert
\right\}
+ \\
\endaligned
\end{equation}
$$
\aligned
& \
\ \ \ \ \
+
y_{x^1}\,y_{x^1}
\cdot
\left\{
-
\left\vert
\begin{array}{ccc}
X_{x^2y}^1 & X_{x^2}^1 & X_y^1 \\
X_{x^2y}^2 & X_{x^2}^2 & X_y^2 \\
Y_{x^2y} & Y_{x^2} & Y_y \\
\end{array}
\right\vert
\right\}
+ \\
& \
\ \ \ \ \ 
+
y_{x^1}\,y_{x^2}
\cdot
\left\{
\left\vert
\begin{array}{ccc}
X_{x^1}^1 & X_{x^2}^1 & X_{yy}^1 \\
X_{x^1}^2 & X_{x^2}^2 & X_{yy}^2 \\
Y_{x^1} & Y_{x^2} & Y_{yy} \\
\end{array}
\right\vert
-
\left\vert
\begin{array}{ccc}
X_{x^1y}^1 & X_{x^2}^1 & X_y^1 \\
X_{x^1y}^2 & X_{x^2}^2 & X_y^2 \\
Y_{x^1y} & Y_{x^2} & Y_y \\
\end{array}
\right\vert
-
\right. \\
& \
\ \ \ \ \
\left.
-
\left\vert
\begin{array}{ccc}
X_{x^1}^1 & X_{x^2y}^1 & X_y^1 \\
X_{x^1}^2 & X_{x^2y}^2 & X_y^2 \\
Y_{x^1} & Y_{x^2y} & Y_y \\
\end{array}
\right\vert
\right\}
+
y_{x^2}\,y_{x^2}
\left\{
-
\left\vert
\begin{array}{ccc}
X_{x^1}^1 & X_{x^1y}^1 & X_y^1 \\
X_{x^1}^2 & X_{x^1y}^2 & X_y^2 \\
Y_{x^1} & Y_{x^1y} & Y_y \\
\end{array}
\right\vert
\right\}
+ \\
& \
\ \ \ \ \
+
y_{x^1}\,y_{x^1}\,y_{x^2}
\cdot
\left\{
-
\left\vert
\begin{array}{ccc}
X_{yy}^1 & X_{x^2}^1 & X_y^1 \\
X_{yy}^2 & X_{x^2}^2 & X_y^2 \\
Y_{yy} & Y_{x^2} & Y_y \\
\end{array}
\right\vert
\right\}
+ 
\\
& \
\ \ \ \ \
+
y_{x^1}\,y_{x^2}\,y_{x^2}
\cdot
\left\{
-
\left\vert
\begin{array}{ccc}
X_{x^1}^1 & X_{yy}^1 & X_y^1 \\
X_{x^1}^2 & X_{yy}^2 & X_y^2 \\
Y_{x^1} & Y_{yy} & Y_y \\
\end{array}
\right\vert
\right\}.
\endaligned
$$ 
Using the third equation~\thetag{ 2.15}, we also deduce, {\it without
any further computation}, that the third equation which would be
obtained by developing~\thetag{ 2.9} in length, is:
\def\theequation{2.23}\begin{equation}
\aligned
0
& 
= 
y_{x^2x^2}
\cdot
\left\vert
\begin{array}{ccc}
X_{x^1}^1 & X_{x^2}^1 & X_y^1 \\
X_{x^1}^2 & X_{x^2}^2 & X_y^2 \\
Y_{x^1} & Y_{x^2} & Y_y \\
\end{array}
\right\vert
+
\left\vert
\begin{array}{ccc}
X_{x^1}^1 & X_{x^2}^1 & X_{x^2x^2}^1 \\
X_{x^1}^2 & X_{x^2}^2 & X_{x^2x^2}^2 \\
Y_{x^1} & Y_{x^2} & Y_{x^2x^2} \\
\end{array}
\right\vert
+ \\
& \
\ \ \ \ \
+
y_{x^1}
\cdot
\left\{
-
\left\vert
\begin{array}{ccc}
X_{x^2x^2}^1 & X_{x^2}^1 & X_y^1 \\
X_{x^1x^2}^2 & X_{x^2}^2 & X_y^2 \\
Y_{x^2x^2} & Y_{x^2} & Y_y \\
\end{array}
\right\vert
\right\}
+ \\
& \
\ \ \ \ \
+
y_{x^2}
\cdot
\left\{
2\,
\left\vert
\begin{array}{ccc}
X_{x^1}^1 & X_{x^2}^1 & X_{x^2y}^1 \\
X_{x^1}^2 & X_{x^2}^2 & X_{x^2y}^2 \\
Y_{x^1} & Y_{x^2} & Y_{x^2y} \\
\end{array}
\right\vert
-
\left\vert
\begin{array}{ccc}
X_{x^1}^1 & X_{x^2x^2}^1 & X_y^1 \\
X_{x^1}^2 & X_{x^2x^2}^2 & X_y^2 \\
Y_{x^1} & Y_{x^2x^2} & Y_y \\
\end{array}
\right\vert
\right\}
+ \\
\endaligned
\end{equation}
$$
\aligned
& \
\ \ \ \ \ 
+
y_{x^1}\,y_{x^2}
\cdot
\left\{
-2\,
\left\vert
\begin{array}{ccc}
X_{x^2y}^1 & X_{x^2}^1 & X_y^1 \\
X_{x^2y}^2 & X_{x^2}^2 & X_y^2 \\
Y_{x^2y} & Y_{x^2} & Y_y \\
\end{array}
\right\vert
\right\}
+ \\
& \
\ \ \ \ \
+
y_{x^2}\,y_{x^2}
\cdot
\left\{
\left\vert
\begin{array}{ccc}
X_{x^1}^1 & X_{x^2}^1 & X_{yy}^1 \\
X_{x^1}^2 & X_{x^2}^2 & X_{yy}^2 \\
Y_{x^1} & Y_{x^2} & Y_{yy} \\
\end{array}
\right\vert
-
2\,
\left\vert
\begin{array}{ccc}
X_{x^1}^1 & X_{x^2y}^1 & X_y^1 \\
X_{x^1}^2 & X_{x^2y}^2 & X_y^2 \\
Y_{x^1} & Y_{x^2y} & Y_y \\
\end{array}
\right\vert
\right\}
+ \\
& \
\ \ \ \ \
+
y_{x^1}\,y_{x^2}\,y_{x^2}
\cdot
\left\{
-
\left\vert
\begin{array}{ccc}
X_{yy}^1 & X_{x^2}^1 & X_y^1 \\
X_{yy}^2 & X_{x^2}^2 & X_y^2 \\
Y_{yy} & Y_{x^2} & Y_y \\
\end{array}
\right\vert
\right\}
+ \\
& \
\ \ \ \ \
+
y_{x^2}\,y_{x^2}\,y_{x^2}
\cdot
\left\{
-
\left\vert
\begin{array}{ccc}
X_{x^1}^1 & X_{yy}^1 & X_y^1 \\
X_{x^1}^2 & X_{yy}^2 & X_y^2 \\
Y_{x^1} & Y_{yy} & Y_y \\
\end{array}
\right\vert
\right\}.
\endaligned
$$

\subsection*{ 2.24.~Appropriate formalism}
To describe the combinatorics underlying
formulas~\thetag{ 2.10}, \thetag{ 2.22} and~\thetag{ 2.23}, 
as in~\cite{ m2004a}, 
let us introduce the following notation 
for the Jacobian determinant:
\def\theequation{2.25}\begin{equation}
\Delta(x^1\vert x^2\vert y)
:=
\left\vert
\begin{array}{ccc}
X_{x^1}^1 & X_{x^2}^1 & X_y^1 \\
X_{x^1}^2 & X_{x^2}^2 & X_y^2 \\
Y_{x^1} & Y_{x^2} & Y_y
\end{array}
\right\vert.
\end{equation}
Here, in the notation $\Delta(x^1\vert x^2\vert y)$, the three spaces
between the two vertical lines $\vert$ refer to the three columns of
the Jacobian determinant, and the terms $x^1$, $x^2$, $y$ in $(x^1
\vert x^2 \vert y)$ designate the partial derivatives appearing in
each column. Accordingly, in the following two examples of {\sl modified
Jacobian determinants}:
\def\theequation{2.26}\begin{equation}
\left\{
\aligned
&
\Delta( \underline{ x^1 x^2} \vert x^2 \vert y)
:=
\left\vert
\begin{array}{ccc}
X_{\underline{x^1x^2}}^1 & X_{x^2}^1 & X_y^1 \\
X_{\underline{x^1x^2}}^2 & X_{x^2}^2 & X_y^2 \\
Y_{\underline{x^1x^2}} & Y_{x^2} & Y_y
\end{array}
\right\vert 
\ \ \ \ \
\ \ \ \ \
{\rm and}
\\
&
\Delta(x^1\vert x^2\vert \underline{x^1y})
:=
\left\vert
\begin{array}{ccc}
X_{x^1}^1 & X_{x^2}^1 & X_{\underline{x^1y}}^1 \\
X_{x^1}^2 & X_{x^2}^2 & X_{\underline{x^1y}}^2 \\
Y_{x^1} & Y_{x^2} & Y_{\underline{x^1y}}
\end{array}
\right\vert,
\endaligned\right.
\end{equation}
we simply mean which column of first order derivatives
is replaced by a column of second order derivatives
in the original Jacobian determinant.

As there are $6$ possible second order derivatives $(\cdot )_{ x^1
x^1}$, $(\cdot )_{ x^1 x^2}$, $(\cdot )_{ x^1 x^y}$, $(\cdot )_{ x^2
x^2}$, $(\cdot )_{ x^2 y}$ and $(\cdot )_{ y y}$ together with $3$
columns, we obtain $3\times 6 = 18$ possible modified Jacobian
determinants:
\def\theequation{2.27}\begin{equation}
\left\{
\aligned
& \
\Delta(x^1x^1 \vert x^2 \vert y) \ \ \ \ \
&
\Delta(x^1 \vert x^1x^1 \vert y) \ \ \ \ \ 
&
\Delta(x^1 \vert x^2 \vert x^1x^1) \ \ \ \ \ 
\\
& \
\Delta(x^1x^2 \vert x^2 \vert y) \ \ \ \ \
&
\Delta(x^1 \vert x^1x^2 \vert y) \ \ \ \ \ 
&
\Delta(x^1 \vert x^2 \vert x^1x^2) \ \ \ \ \ 
\\
& \
\Delta(x^1y \vert x^2 \vert y) \ \ \ \ \
&
\Delta(x^1 \vert x^1y \vert y) \ \ \ \ \ 
&
\Delta(x^1 \vert x^2 \vert x^1y) \ \ \ \ \ 
\\
& \
\Delta(x^2x^2 \vert x^2 \vert y) \ \ \ \ \
&
\Delta(x^1 \vert x^2x^2 \vert y) \ \ \ \ \ 
&
\Delta(x^1 \vert x^2 \vert x^2x^2) \ \ \ \ \ 
\\
& \
\Delta(x^2y \vert x^2 \vert y) \ \ \ \ \
&
\Delta(x^1 \vert x^2y \vert y) \ \ \ \ \ 
&
\Delta(x^1 \vert x^2 \vert x^2y) \ \ \ \ \ 
\\
& \
\Delta(yy \vert x^2 \vert y) \ \ \ \ \
&
\Delta(x^1 \vert yy \vert y) \ \ \ \ \ 
&
\Delta(x^1 \vert x^2 \vert yy). \ \ \ \ \ 
\\
\endaligned\right.
\end{equation}

Next, we observe that if we want to solve with
respect to $y_{ x^1x^1}$ in~\thetag{ 2.10}, 
with respect to $y_{ x^1x^2}$ in~\thetag{ 2.22}
and with respect to $y_{ x^2x^2}$ in~\thetag{ 2.23}, 
we have to divide by the Jacobian determinant 
$\Delta( x^1 \vert x^2\vert y)$. 
Consequently, we introduce 
$18$ new {\sl square functions} as follows:
\def\theequation{2.28}\begin{equation}
\small
\left\{
\aligned
\square_{x^1x^1}^1:=&\ \frac{ \Delta(x^1x^1\vert x^2 \vert y)}{
\Delta(x^1\vert x^2 \vert y)} 
\ \ \ \ \ & 
\square_{x^1x^2}^1:=&\ \frac{ \Delta(x^1x^2\vert x^2 \vert y)}{
\Delta(x^1\vert x^2 \vert y)} 
\ \ \ \ \ & 
\square_{x^1y}^1:=&\ \frac{ \Delta(x^1y\vert x^2 \vert y)}{
\Delta(x^1\vert x^2 \vert y)} 
\ \ \ \ \ 
\\
\square_{x^2x^2}^1:=&\ \frac{ \Delta(x^2x^2\vert x^2 \vert y)}{
\Delta(x^1\vert x^2 \vert y)} 
\ \ \ \ \ & 
\square_{x^2y}^1:=&\ \frac{ \Delta(x^2y\vert x^2 \vert y)}{
\Delta(x^1\vert x^2 \vert y)} 
\ \ \ \ \ & 
\square_{yy}^1:=&\ \frac{ \Delta(yy\vert x^2 \vert y)}{
\Delta(x^1\vert x^2 \vert y)} 
\ \ \ \ \ 
\\
\square_{x^1x^1}^2:=&\ \frac{ \Delta(x^1\vert x^1x^1 \vert y)}{
\Delta(x^1\vert x^2 \vert y)} 
\ \ \ \ \ & 
\square_{x^1x^2}^2:=&\ \frac{ \Delta(x^1\vert x^1x^2 \vert y)}{
\Delta(x^1\vert x^2 \vert y)} 
\ \ \ \ \ & 
\square_{x^1y}^2:=&\ \frac{ \Delta(x^1\vert x^1y \vert y)}{
\Delta(x^1\vert x^2 \vert y)} 
\ \ \ \ \
\\
\square_{x^2x^2}^2:=&\ \frac{ \Delta(x^1\vert x^2x^2 \vert y)}{
\Delta(x^1\vert x^2 \vert y)} 
\ \ \ \ \ & 
\square_{x^2y}^2:=&\ \frac{ \Delta(x^1\vert x^2y \vert y)}{
\Delta(x^1\vert x^2 \vert y)} 
\ \ \ \ \ & 
\square_{yy}^2:=&\ \frac{ \Delta(x^1\vert yy \vert y)}{
\Delta(x^1\vert x^2 \vert y)} 
\ \ \ \ \ 
\\
\square_{x^1x^1}^3:=&\ \frac{ \Delta(x^1\vert x^2 \vert x^1x^1)}{
\Delta(x^1\vert x^2 \vert y)} 
\ \ \ \ \ & 
\square_{x^1x^2}^3:=&\ \frac{ \Delta(x^1\vert x^2 \vert x^1x^2)}{
\Delta(x^1\vert x^2 \vert y)} 
\ \ \ \ \ & 
\square_{x^1y}^3:=&\ \frac{ \Delta(x^1\vert x^2 \vert x^1y)}{
\Delta(x^1\vert x^2 \vert y)} 
\ \ \ \ \ 
\\
\square_{x^2x^2}^3:=&\ \frac{ \Delta(x^1\vert x^2 \vert x^2x^2)}{
\Delta(x^1\vert x^2 \vert y)} 
\ \ \ \ \ & 
\square_{x^2y}^3:=&\ \frac{ \Delta(x^1\vert x^2 \vert x^2y)}{
\Delta(x^1\vert x^2 \vert y)} 
\ \ \ \ \ & 
\square_{yy}^3:=&\ \frac{ \Delta(x^1\vert x^2 \vert yy)}{
\Delta(x^1\vert x^2 \vert y)}.
\ \ \ \ \ 
\\
\endaligned\right.
\end{equation} 

Thanks to these notations, we can rewrite the three equations~\thetag{
2.10}, \thetag{ 2.22} and~\thetag{ 2.23} in a more compact style.

\def\thelemma{2.29}\begin{lemma}
A completely integrable system of
{\rm three} second order partial 
differential equations 
\def\theequation{2.30}\begin{equation}
\left\{
\aligned
y_{x^1x^1} (x)
&
= 
F^{1,1}
\left(
x^1, x^2, y(x), y_{x^1}(x), y_{x^2}(x) 
\right), 
\\
y_{x^1x^2} (x)
& 
= 
F^{1,2}
\left(
x^1, x^2, y(x), y_{x^1}(x), y_{x^2}(x) 
\right), 
\\
y_{x^2x^2} (x)
&
= 
F^{2,2}
\left(
x^1, x^2, y(x), y_{x^1}(x), y_{x^2}(x) 
\right), 
\\
\endaligned\right.
\end{equation}
is equivalent to the simplest
system $Y_{X^1 X^1} = 0$, 
$Y_{ X^1 X^2} = 0$, $Y_{X^2, X^2} = 0$,
{\rm if and only if}
there exist local 
$\K$-analytic functions 
$X^1$, $X^2$, $Y$ such that it
may be written under the
specific form{\rm :}
\def\theequation{2.31}\begin{equation}
\small
\left\{
\aligned
y_{x^1x^1}
& 
=
- 
\square_{x^1x^1}^3
+
y_{x^1}
\cdot
\left(
-
2\,\square_{x^1y}^3
+
\square_{x^1x^1}^1
\right)
+
y_{x^2}
\cdot
\left(
\square_{x^1x^1}^2
\right)
+ \\
& \
\ \ \ \ \ 
+
y_{x^1}\,y_{x^1}
\cdot
\left(
-
\square_{yy}^3
+
2\,\square_{x^1y}^1
\right)
+
y_{x^1}\,y_{x^2}
\cdot
\left(
2\,\square_{x^1y}^2
\right)
+ \\
& \
\ \ \ \ \ 
+
y_{x^1}\,y_{x^1}\,y_{x^1}
\cdot
\left(
\square_{yy}^1
\right)
+
y_{x^1}\,y_{x^1}\,y_{x^2}
\cdot
\left(
\square_{yy}^2
\right),
\\
y_{x^1x^2}
& 
=
- 
\square_{x^1x^2}^3
+
y_{x^1}
\cdot
\left(
-
\square_{x^2y}^3
+
\square_{x^1x^2}^1
\right)
+
y_{x^2}
\cdot
\left(
-
\square_{x^1y}^3
+
\square_{x^1x^2}^2
\right)
+ \\
& \
\ \ \ \ \
+
y_{x^1}\,y_{x^1}
\cdot
\left(
\square_{x^2y}^1
\right)
+
y_{x^1}\,y_{x^2}
\cdot
\left(
-
\square_{yy}^3
+
\square_{x^1y}^1
+
\square_{x^2y}^2
\right)
+ \\
& \
\ \ \ \ \ 
+
y_{x^2}\,y_{x^2}
\cdot
\left(
\square_{x^1y}^2
\right)
+
y_{x^1}\,y_{x^1}\,y_{x^2}
\cdot
\left(
\square_{yy}^1
\right)
+
y_{x^1}\,y_{x^2}\,y_{x^2}
\cdot
\left(
\square_{yy}^2
\right),
\\
y_{x^2x^2}
& 
=
- 
\square_{x^2x^2}^3
+
y_{x^1}
\cdot
\left(
\square_{x^2x^2}^1
\right)
+
y_{x^2}
\cdot
\left(
-
2\,\square_{x^2y}^3
+
\square_{x^2x^2}^2
\right)
+ \\
& \
\ \ \ \ \ 
+
y_{x^1}\,y_{x^2}
\cdot
\left(
2\,\square_{x^2y}^1
\right)
+
y_{x^2}\,y_{x^2}
\cdot
\left(
-
\square_{yy}^3
+
2\,\square_{x^2y}^2
\right)
+ \\
& \
\ \ \ \ \ 
+
y_{x^1}\,y_{x^2}\,y_{x^2}
\cdot
\left(
\square_{yy}^1
\right)
+
y_{x^2}\,y_{x^2}\,y_{x^2}
\cdot
\left(
\square_{yy}^2
\right).
\endaligned\right.
\end{equation}
\end{lemma}

\subsection*{ 2.32.~General formulas}
The formal dictionary between the original determinantial
formulas~\thetag{ 2.10}, \thetag{ 2.22}, \thetag{ 2.23}, between the
coefficients~\thetag{ 2.15} of the second order prolongation of a
vector field and between the new square formulas~\thetag{ 2.31} above
is evident. Consequently, {\it without any computation}, just by
translating the family of formulas~\thetag{ 2.14}, we may deduce the
exact formulation of the desired generalization of Lemma~2.29 above.

\def\thelemma{2.33}\begin{lemma}
A completely integrable system of second order partial 
differential equations of the form
\def\theequation{2.34}\begin{equation}
y_{x^{j_1}x^{j_2}}(x) 
=
F^{j_1, j_2} 
\left(
x, y(x), y_{x^1}(x),\dots,y_{x^n}(x)
\right), 
\ \ \ \ \ \ \ \ 
j_1,j_2=1,\dots n, 
\end{equation}
is equivalent to the simplest system $Y_{X^{j_1} X^{j_2}} = 0$, $j_1,
j_2 = 1, \dots, n$, {\rm if and only if} there exist local
$\K$-analytic functions $X^{l}$, $Y$ such that it may be written under
the specific form{\rm :}
\def\theequation{2.35}\begin{equation}
\left\{
\aligned
y_{x^{j_1}x^{j_2}}
&
=
-
\square_{x^{j_1}x^{j_2}}^{n+1}
+
\sum_{k_1=1}^n\,y_{x^{k_1}}
\cdot
\left\{
\left(
\square_{x^{j_1}x^{j_2}}^{k_1}
-
\delta_{j_1}^{k_1}\,\square_{x^{j_2}y}^{n+1}
-
\delta_{j_2}^{k_1}\,\square_{x^{j_1}y}^{n+1}
\right)
+ 
\right.
\\
& \
\ \ \ \ \
\left.
+
y_{x^{j_1}}
\cdot
\left(
\square_{x^{j_2}y}^{k_1}
-
\frac{1}{2}\,\delta_{j_2}^{k_1}\,\square_{yy}^{n+1}
\right)
+
y_{x^{j_2}}
\cdot
\left(
\square_{x^{j_1}y}^{k_1}
-
\frac{1}{2}\,\delta_{j_1}^{k_1}\,\square_{yy}^{n+1}
\right)
+ 
\right.
\\
& \
\ \ \ \ \
\left.
+
y_{x^{j_1}}\,y_{x^{j_2}}
\cdot
\left(
\square_{yy}^{k_1}
\right)
\right\}.
\endaligned\right.
\end{equation}
\end{lemma}

Of course, to define the square functions in the context of $n\geq 2$
independent variables $(x^1, x^2, \dots, x^n)$, we introduce the
Jacobian determinant
\def\theequation{2.36}\begin{equation}
\Delta(x^1 \vert x^2 \vert \cdots \vert x^n \vert y)
:= 
\left\vert
\begin{array}{cccc}
X_{x^1}^1 & \cdots & X_{x^n}^1 & X_y^1 \\
\vdots & \cdots & \vdots & \vdots \\
X_{x^1}^n & \cdots & X_{x^n}^n & X_y^n \\
Y_{x^1} & \cdots & Y_{x^n} & Y_y \\
\end{array}
\right\vert,
\end{equation}
together with its modifications
\def\theequation{2.37}\begin{equation}
\Delta
\left(
x^1\vert\cdots\vert^{k_1}\,x^{j_1}\,x^{j_2} \vert \cdots \vert y
\right),
\end{equation}
in which the $k_1$-th column of partial first order derivatives
$\vert^{k_1}\, x^{k_1} \vert$ is replaced by the column $\vert^{k_1}
\, x^{j_1} x^{j_2} \vert$ of partial derivatives. Here, the indices
$k_1$, $j_1$, $j_2$ satisfy $1\leq k_1, j_1, j_2 \leq n+1$, with the
convention that we adopt the notational equivalence
\def\theequation{2.38}\begin{equation}
\fbox{
$x^{n+1}\equiv y$
}.
\end{equation}
This convention will be convenient to write some of our general
formulas in the sequel.

As we promised to only summarize the proof of Theorem~1.1 in this
paper, we will not reproduce the proof of Lemma~2.37 from~\cite{
m2003}. Involving linear algebra considerations, this proof is similar
to (and in fact slightly simpler than) the proof of the analogous
Lemma~3.32 in~\cite{ m2004a}.

\section*{ \S3.~First and second auxiliary system}

\subsection*{3.1.~Functions $G_{j_1, j_2}$, 
$H_{j_1, j_2}^{k_1}$, $L_{j_1}^{k_1}$ and $M^{k_1}$} To discover the
four families of functions appearing in the statement of Theorem~1.1,
by comparing~\thetag{ 2.35} and~\thetag{ 1.3},
it suffices (of course) to set:
\def\theequation{3.2}\begin{equation}
\left\{
\aligned
G_{j_1,j_2}
&
:= 
-
\square_{x^{j_1}x^{j_2}}^{n+1}, 
\\
H_{j_1,j_2}^{k_1}
&
:=
\square_{x^{j_1}x^{j_2}}^{k_1}
-
\delta_{j_1}^{k_1}\,\square_{x^{j_2}y}^{n+1}
-
\delta_{j_2}^{k_1}\,\square_{x^{j_1}y}^{n+1}, 
\\
L_{j_1}^{k_1}
&
:=
2\,\square_{x^{j_1}y}^{k_1}
-
\delta_{j_1}^{k_1}\,\square_{yy}^{n+1},
\\
M^{k_1}
&
:=
\square_{yy}^{k_1}.
\endaligned\right.
\end{equation}

Consequently, we have shown the ``only if'' part of Theorem~1.1, which
is the easiest implication.

To establish the ``if'' part, by far the most difficult implication,
the very main lemma can be stated as follows.

\def\thelemma{3.3}\begin{lemma}
The partial differerential relations {\rm (I')}, {\rm (II')}, {\rm
(III')} and {\rm (IV')} which express in length the compatibility
conditions for the system~\thetag{ 1.3} are necessary and sufficient
for the existence of functions $X^l$, $Y$ of $(x^{l_1}, y)$ satisfying
the second order nonlinear system of partial differential
equations~\thetag{ 3.2} above.
\end{lemma}

Indeed, the collection of equations~\thetag{ 3.2} is a system of
partial differential equations with unknowns $X^l$, $Y$, by virtue of
the definition of the square functions.

\subsection*{ 3.4.~First auxiliary system}
To proceed further, we observe that there are $(m+1)$ more square
functions than functions $G_{j_1, j_2}$, $H_{j_1, j_2}^{k_1}$,
$L_{j_1}^{k_1}$ and $M^{k_1}$. Indeed, a simple counting yields:
\def\theequation{3.5}\begin{equation}
\left\{
\aligned
{}
& \
\# \{\square_{x^{j_1}x^{j_2}}^{k_1}\}
=
\frac{n^2(n+1)}{2},
& \
\# \{\square_{x^{j_1}y}^{k_1} \}
=
n^2, \\
& \
\# \{\square_{yy}^{k_1}\}
=
n, 
& \
\# \{\square_{x^{j_1}x^{j_2}}^{n+1} \}
= 
\frac{n(n+1)}{2}, \\
& \
\# \{\square_{x^{j_1}y}^{n+1}\}
= 
n, 
\ \ \ \ \ \ \ \ \ \ \ 
\ \ \ \ \ \ \ \ \ \ \ 
& \
\# \{\square_{yy}^{n+1} \}
= 
1, \\
\endaligned\right.
\end{equation}
whereas
\def\theequation{3.6}\begin{equation}
\left\{
\aligned
{}
& \
\# \{G^{j_1,j_2}\}
= 
\frac{n(n+1)}{2}, 
& \ 
\#
\{H_{j_1,j_2}^{k_1} \}
=
\frac{n^2(n+1)}{2}, \\
& \
\# \{L_{j_1}^{k_1} \}
=
n^2, 
& \ 
\ \ \ \ \ \ \ \ \ \ \ 
\ \ \ \ \ \ \ \ \ \ \
\# \{M^{k_1}\}
=
n.
\endaligned\right.
\end{equation}
Here, the indices $j_1$, $j_2$, $k_1$ satisfy $1\leq j_1, j_2, k_1
\leq n$. Similarly as in~\cite{ m2004a}, to transform the
system~\thetag{ 3.2} in a true complete system, let us introduce
functions $\Pi_{j_1,j_2}^{k_1}$ of $(x^{l_1}, y)$, where $1\leq j_1,
j_2, k_1 \leq n+1$, which satisfy the symmetry $\Pi_{j_1,j_2}^{k_1} =
\Pi_{j_1, j_1}^{k_1}$, and let us introduce the
following {\sl first auxiliary
system}:
\def\theequation{3.7}\begin{equation}
\left\{
\aligned
\square_{x^{j_1}x^{j_2}}^{k_1}
=
\Pi_{j_1,j_2}^{k_1}, 
\ \ \ \ \ \ \ &
\square_{x^{j_1}y}^{k_1}
=
\Pi_{j_1, n+1}^{k_1}, 
\ \ \ \ \ \ \ &
\square_{yy}^{k_1}
=
\Pi_{n+1,n+1}^{k_1}, 
\\
\square_{x^{j_1}x^{j_2}}^{n+1}
=
\Pi_{j_1,j_2}^{n+1}, 
\ \ \ \ \ \ \ &
\square_{x^{j_1}y}^{n+1}
=
\Pi_{j_1, n+1}^{n+1}, 
\ \ \ \ \ \ \ &
\square_{yy}^{n+1}
=
\Pi_{n+1,n+1}^{n+1}. 
\endaligned\right.
\end{equation}
It is complete. The necessary and sufficient
conditions for the existence of solutions
$X^l$, $Y$ follow by cross
differentiations.

\def\thelemma{3.8}\begin{lemma}
For all $j_1, j_2, j_3, k_1 = 1, 2, \dots, n+1$, 
we have the cross differentiation relations
\def\theequation{3.9}\begin{equation}
\left(
\square_{x^{j_1}x^{j_2}}^{k_1}
\right)_{x^{j_3}}
-
\left(
\square_{x^{j_1}x^{j_3}}^{k_1}
\right)_{x^{j_2}}
=
-
\sum_{k_2=1}^{n+1}\,
\square_{x^{j_1}x^{j_2}}^{k_2}\,\square_{x^{j_3}x^{k_2}}^{k_1}
+
\sum_{k_2=1}^{n+1}\,
\square_{x^{j_1}x^{j_3}}^{k_2}\,\square_{x^{j_2}x^{k_2}}^{k_1}.
\end{equation}
\end{lemma}

The proof of this lemma is exactly the same as the proof of Lemma~3.40
in~\cite{ m2004a}.

As a direct consequence, we deduce that a necessary and sufficient
condition for the existence of solutions $\Pi_{j_1, j_2}^{ k_1}$ to
the first auxiliary system is that they satisfy the following
compatibility partial differential relations:
\def\theequation{3.10}\begin{equation}
\frac{\partial \Pi_{j_1,j_2}^{k_1}}{\partial x^{j_3}}
-
\frac{\partial \Pi_{j_1,j_3}^{k_1}}{\partial x^{j_2}}
=
-
\sum_{k_2=1}^{n=1}\,\Pi_{j_1,j_2}^{k_2}
\cdot
\Pi_{j_3,k_2}^{k_1}
+
\sum_{k_2=1}^{n=1}\,\Pi_{j_1,j_3}^{k_2}
\cdot
\Pi_{j_2,k_2}^{k_1},
\end{equation}
for all $j_1, j_2, j_3, k_1 = 1, \dots, n+1$.

We shall have to specify this system in length according to the
splitting $\{1, 2, \dots, n\}$ and $\{n+1\}$ of the indices of
coordinates. We obtain six families of equations equivalent
to~\thetag{ 3.10} just above:
\def\theequation{3.11}\begin{equation}
\left\{
\aligned
\left(
\Pi_{j_1,j_2}^{n+1}
\right)_{x^{j_3}}
-
\left(
\Pi_{j_1,j_3}^{n+1}
\right)_{x^{j_2}}
&
=
-
\sum_{k_2=1}^n\,\Pi_{j_1,j_2}^{k_2}\,\Pi_{j_3,k_2}^{n+1}
-
\Pi_{j_1,j_2}^{n+1}\,\Pi_{j_3,n+1}^{n+1}
+ \\
& \
\ \ \ \ \ \ \ \ \ \
+
\sum_{k_2=1}^n\,\Pi_{j_1,j_3}^{k_2}\,\Pi_{j_2,k_2}^{n+1}
+
\Pi_{j_1,j_3}^{n+1}\,\Pi_{j_2,n+1}^{n+1}, 
\\
\left(
\Pi_{j_1,j_2}^{n+1}
\right)_y
-
\left(
\Pi_{j_1,n+1}^{n+1}
\right)_{x^{j_2}}
&
=
-
\sum_{k_2=1}^n\,\Pi_{j_1,j_2}^{k_2}\,\Pi_{n+1,k_2}^{n+1}
-
\Pi_{j_1,j_2}^{n+1}\,\Pi_{n+1,n+1}^{n+1}
+ \\
& \
\ \ \ \ \ \ \ \ \ \
+
\sum_{k_2=1}^n\,\Pi_{j_1,n+1}^{k_2}\,\Pi_{j_2,k_2}^{n+1}
+
\Pi_{j_1,n+1}^{n+1}\,\Pi_{j_2,n+1}^{n+1}, 
\\
\left(
\Pi_{j_1,n+1}^{n+1}
\right)_y
-
\left(
\Pi_{n+1,n+1}^{n+1}
\right)_{x^{j_1}}
&
=
-
\sum_{k_2=1}^n\,\Pi_{j_1,n+1}^{k_2}\,\Pi_{n+1,k_2}^{n+1}
-
\underline{
\Pi_{j_1,n+1}^{n+1}\,\Pi_{n+1,n+1}^{n+1} 
}_{ \octagon \! \! \! \! \! \tiny{\sf a}}
+ \\
& \
\ \ \ \ \ \ \ \ \ \
+
\sum_{k_2=1}^n\,\Pi_{n+1,n+1}^{k_2}\,\Pi_{j_1,k_2}^{n+1}
+
\underline{
\Pi_{n+1,n+1}^{n+1}\,\Pi_{j_1,n+1}^{n+1}
}_{ \octagon \! \! \! \! \! \tiny{\sf a}} \ , 
\\
\left(
\Pi_{j_1,j_2}^{k_1}
\right)_{x^{j_3}}
-
\left(
\Pi_{j_1,j_3}^{k_1}
\right)_{x^{j_2}}
&
=
-
\sum_{k_2=1}^n\,\Pi_{j_1,j_2}^{k_2}\,\Pi_{j_3,k_2}^{k_1}
-
\Pi_{j_1,j_2}^{n+1}\,\Pi_{j_3,n+1}^{k_1}
+ \\
& \
\ \ \ \ \ \ \ \ \ \
+
\sum_{k_2=1}^n\,\Pi_{j_1,j_3}^{k_2}\,\Pi_{j_2,k_2}^{k_1}
+
\Pi_{j_1,j_3}^{n+1}\,\Pi_{j_2,n+1}^{k_1}, 
\\
\left(
\Pi_{j_1,j_2}^{k_1}
\right)_y
-
\left(
\Pi_{j_1,n+1}^{k_1}
\right)_{x^{j_2}}
&
=
-
\sum_{k_2=1}^n\,\Pi_{j_1,j_2}^{k_2}\,\Pi_{n+1,k_2}^{k_1}
-
\Pi_{j_1,j_2}^{n+1}\,\Pi_{n+1,n+1}^{k_1}
+ \\
& \
\ \ \ \ \ \ \ \ \ \
+
\sum_{k_2=1}^n\,\Pi_{j_1,n+1}^{k_2}\,\Pi_{j_2,k_2}^{k_1}
+
\Pi_{j_1,n+1}^{n+1}\,\Pi_{j_2,n+1}^{k_1}, 
\\
\left(
\Pi_{j_1,n+1}^{k_1}
\right)_y
-
\left(
\Pi_{n+1,n+1}^{k_1}
\right)_{x^{j_1}}
&
=
-
\sum_{k_2=1}^n\,\Pi_{j_1,n+1}^{k_2}\,\Pi_{n+1,k_2}^{k_1}
-
\Pi_{j_1,n+1}^{n+1}\,\Pi_{n+1,n+1}^{k_1}
+ \\
& \
\ \ \ \ \ \ \ \ \ \
+
\sum_{k_2=1}^n\,\Pi_{n+1,n+1}^{k_2}\,\Pi_{j_1,k_2}^{k_1}
+
\Pi_{n+1,n+1}^{n+1}\,\Pi_{j_1,n+1}^{k_1}. 
\\
\endaligned\right.
\end{equation}
where the indices $j_1, j_2, j_3, k_1$ vary in the set $\{ 1, 2, 1,
\dots, n\}$.

\subsection*{3.12.~Principal unknowns}
As there are $(m+1)$ more square (or Pi) functions than the functions
$G_{j_1, j_2}$, $H_{j_1, j_2}^{k_1}$, $L_{j_1}^{k_1}$ and $M^{ k_1}$,
we cannot invert directly the linear system~\thetag{ 3.2}. To
quasi-inverse it, we choose the $(m+1)$ specific square functions
\def\theequation{3.13}\begin{equation}
\Theta^1
:= 
\square_{x^1x^1}^1, 
\ \ \ \ \
\Theta^2
:= 
\square_{x^2x^2}^2,
\cdots\cdots,
\Theta^{n+1}
:= 
\square_{x^{n+1}x^{n+1}}^{n+1},
\end{equation}
calling them {\sl principal unknowns}, 
and we get the quasi-inversion:
\def\theequation{3.14}\begin{equation}
\left\{
\aligned
\Pi_{j_1,j_2}^{k_1}
&
=
\square_{x^{j_1}x^{j_2}}^{k_1}
=
H_{j_1,j_2}^{k_1}
-
\frac{1}{2}\,\delta_{j_1}^{k_1}\,H_{j_2,j_2}^{j_2}
-
\frac{1}{2}\,\delta_{j_2}^{k_1}\,H_{j_1,j_1}^{j_1}
+
\frac{1}{2}\,\delta_{j_1}^{k_1}\,\Theta^{j_2}
+
\frac{1}{2}\,\delta_{j_2}^{k_1}\,\Theta^{j_1}, 
\\
\Pi_{j_1,n+1}^{k_1}
&
=
\square_{x^{j_1}y}^{k_1}
=
\frac{1}{2}\,L_{j_1}^{k_1}
+
\frac{1}{2}\,\delta_{j_1}^{k_1}\,\Theta^{n+1}, 
\\
\Pi_{n+1,n+1}^{k_1}
&
=
\square_{yy}^{k_1}
=
M^{k_1},
\\
\Pi_{j_1,j_2}^{n+1}
&
=
\square_{x^{j_1}x^{j_2}}^{n+1}
=
-
G_{j_1,j_2},
\\
\Pi_{j_1,n+1}^{n+1}
&
=
\square_{x^{j_1}y}^{n+1}
=
-
\frac{1}{2}\,H_{j_1, j_1}^{j_1}
+
\frac{1}{2}\,\Theta^{j_1}.
\endaligned\right.
\end{equation}

\subsection*{ 3.15.~Second auxiliary system}
Replacing the five families of functions $\Pi_{j_1, j_2}^{k_1}$,
$\Pi_{j_1, n+1}^{k_1}$, $\Pi_{ n+1, n+1}^{k_1}$, $\Pi_{j_1,
j_2}^{n+1}$, $\Pi_{j_1, n+1}^{n+1}$ by their values obtained
in~\thetag{ 3.14} just above together with the
principal unknowns
\def\theequation{3.16}\begin{equation}
\Pi_{j_1, j_1}^{j_1} 
=
\Theta^{j_1}, 
\ \ \ \ \ \ \ \ \ \
\Pi_{n+1,n+1}^{n+1}
= 
\Theta^{n+1},
\end{equation}
in the six equations $(3.11)_1$, $(3.11)_2$, $(3.11)_3$, $(3.11)_4$,
$(3.11)_5$ and $(3.11)_6$, after hard computations that we will not
reproduce here, we obtain six families of equations. From now on, we
abbreviate every sum $\sum_{k=1}^n$ as $\sum_{k_1}$.

Firstly:
\def\theequation{3.17}\begin{equation}
0
=
\underline{
G_{j_1,j_2,x^{j_3}}
-
G_{j_1,j_3,x^{j_2}}
}
+
\sum_{k_1}\,G_{j_3,k_1}\,H_{j_1,j_2}^{k_1}
-
\sum_{k_1}\,G_{j_2,k_1}\,H_{j_1,j_3}^{k_1}.
\end{equation}
This is (I') of Theorem~1.1.
Just above and below, we plainly underline the
monomials involving a first order derivative. Secondly:
\def\theequation{3.18}\begin{equation}
\left\{
\aligned
\Theta_{x^{j_2}}^{j_1}
&
=
-
\underline{
2\,G_{j_1,j_2,y}
+
H_{j_1,j_1,x^{j_2}}^{j_1}
}
+ \\
& \
\ \ \ \ \
+
\sum_{k_1}\,G_{j_2,k_1}\,L_{j_1}^{k_1}
+
\frac{1}{2}\,H_{j_1,j_1}^{j_1}\,H_{j_2,j_2}^{j_2}
-
\sum_{k_1}\,H_{j_1,j_2}^{k_1}\,H_{k_1,k_1}^{k_1}
- \\
& \
\ \ \ \ \
-
G_{j_1,j_2}\,\Theta^{n+1}
-
\frac{1}{2}\,H_{j_1,j_1}^{j_1}\,\Theta^{j_2}
-
\frac{1}{2}\,H_{j_2,j_2}^{j_2}\,\Theta^{j_1}
+
\sum_{k_1}\,H_{j_1,j_2}^{k_1}\,\Theta^{k_1}
+ \\
& \
\ \ \ \ \
+ 
\frac{1}{2}\,\Theta^{j_1}\,\Theta^{j_2}.
\endaligned\right.
\end{equation}
Thirdly:
\def\theequation{3.19}\begin{equation}
\left\{
\aligned
-
\Theta_{x^{j_1}}^{n+1}
+
\frac{1}{2}\,\Theta_y^{j_1}
& 
=
\underline{
\frac{1}{2}\,H_{j_1,j_1,y}^{j_1}
}
- \\
& \
\ \ \ \ \ \
-
\sum_{k_1}\,G_{j_1,k_1}\,M^{k_1}
+
\frac{1}{4}\,\sum_{k_1}\,H_{k_1,k_1}^{k_1}\,L_{j_1}^{k_1}
+ \\
& \
\ \ \ \ \ 
+
\frac{1}{4}\,H_{j_1,j_1}^{j_1}\,\Theta^{n+1}
-
\frac{1}{4}\,\sum_{k_1}\,L_{j_1}^{k_1}\,\Theta^{k_1}
-
\frac{1}{4}\,\Theta^{j_1}\,\Theta^{n+1}.
\endaligned\right.
\end{equation}
Fourtly:
\def\theequation{3.20}\begin{equation}
\left\{
\aligned
&
\frac{1}{2}\,\delta_{j_1}^{k_1}\,\Theta_{x^{j_3}}^{j_2}
-
\frac{1}{2}\,\delta_{j_1}^{k_1}\,\Theta_{x^{j_2}}^{j_3}
+
\frac{1}{2}\,\delta_{j_2}^{k_1}\,\Theta_{x^{j_3}}^{j_1}
-
\frac{1}{2}\,\delta_{j_3}^{k_1}\,\Theta_{x^{j_2}}^{j_1}
= \\
& \
\ \ \ \ \
=
-
\underline{
H_{j_1,j_2,x^{j_3}}^{k_1}
+
H_{j_1,j_3,x^{j_2}}^{k_1}
-
\frac{1}{2}\,\delta_{j_1}^{k_1}\,H_{j_3,j_3,x^{j_2}}^{j_3}
+
\frac{1}{2}\,\delta_{j_1}^{k_1}\,H_{j_2,j_2,x^{j_3}}^{j_2}
}
- \\
& \
\ \ \ \ \ \ \ \ \ \ 
-
\underline{
\frac{1}{2}\,\delta_{j_3}^{k_1}\,H_{j_1,j_1,x^{j_2}}^{j_1}
+
\frac{1}{2}\,\delta_{j_2}^{k_1}\,H_{j_1,j_1,x^{j_3}}^{j_1}
}
- \\
& \
\ \ \ \ \ \ \ \ \ \ 
-
\frac{1}{2}\,G_{j_1,j_2}\,L_{j_3}^{k_1}
+
\frac{1}{2}\,G_{j_1,j_3}\,L_{j_2}^{k_1}
-
\frac{1}{4}\,\delta_{j_3}^{k_1}\,H_{j_1,j_1}^{j_1}\,H_{j_2,j_2}^{j_2}
+
\frac{1}{4}\,\delta_{j_2}^{k_1}\,H_{j_1,j_1}^{j_1}\,H_{j_3,j_3}^{j_3}
- \\
& \
\ \ \ \ \ \ \ \ \ \ 
-
\sum_{k_2}\,H_{j_1,j_2}^{k_2}\,H_{j_3,k_2}^{k_1}
+
\sum_{k_2}\,H_{j_1,j_3}^{k_2}\,H_{j_2,k_2}^{k_1}
-
\frac{1}{2}\,\delta_{j_2}^{k_1}\,H_{j_1,j_3}^{k_2}\,H_{k_2,k_2}^{k_2}
+
\frac{1}{2}\,\delta_{j_3}^{k_1}\,H_{j_1,j_2}^{k_2}\,H_{k_2,k_2}^{k_2}
- \\
& \
\ \ \ \ \ \ \ \ \ \ 
-
\frac{1}{2}\,\delta_{j_2}^{k_1}\,G_{j_1,j_3}\,\Theta^{n+1}
+
\frac{1}{2}\,\delta_{j_3}^{k_1}\,G_{j_1,j_2}\,\Theta^{n+1}
- \\
& \
\ \ \ \ \ \ \ \ \ \ 
-
\frac{1}{4}\,\delta_{j_2}^{k_1}\,H_{j_1,j_1}^{j_1}\,\Theta^{j_3}
+
\frac{1}{4}\,\delta_{j_3}^{k_1}\,H_{j_1,j_1}^{j_1}\,\Theta^{j_2}
-
\frac{1}{4}\,\delta_{j_2}^{k_1}\,H_{j_3,j_3}^{j_3}\,\Theta^{j_1}
+
\frac{1}{4}\,\delta_{j_3}^{k_1}\,H_{j_2,j_2}^{j_2}\,\Theta^{j_1}
- \\
& \
\ \ \ \ \ \ \ \ \ \ 
-
\frac{1}{2}\,\delta_{j_3}^{k_1}\,
\sum_{k_1}\,H_{j_1,j_2}^{k_2}\,\Theta^{k_2}
+
\frac{1}{2}\,\delta_{j_2}^{k_1}\,
\sum_{k_1}\,H_{j_1,j_3}^{k_2}\,\Theta^{k_2}- \\
& \
\ \ \ \ \ \ \ \ \ \ 
-
\frac{1}{4}\,\delta_{j_3}^{k_1}\,\Theta^{j_1}\,\Theta^{j_2}
+
\frac{1}{4}\,\delta_{j_2}^{k_1}\,\Theta^{j_1}\,\Theta^{j_3}.
\endaligned\right.
\end{equation}
Fifthly:
\def\theequation{3.21}\begin{equation}
\left\{
\aligned
&
\frac{1}{2}\,\delta_{j_1}^{k_1}\,\Theta_y^{j_2}
+
\frac{1}{2}\,\delta_{j_2}^{k_1}\,\Theta_y^{j_1}
-
\frac{1}{2}\,\delta_{j_1}^{k_1}\,\Theta_{x^{j_2}}^{n+1}
= \\
& \
\ \ \ \ \
=
G_{j_1,j_2}\,M^{k_1}
+
\frac{1}{2}\,\sum_{k_2}\,H_{j_1,k_2}^{k_1}\,L_{j_1}^{k_2}
-
\frac{1}{2}\,\sum_{k_2}\,H_{j_1,j_2}^{k_2}\,L_{k_2}^{k_1}
-
\frac{1}{4}\,\delta_{j_2}^{k_1}\,\sum_{k_2}\,
H_{k_2,k_2}^{k_2}\,L_{j_1}^{k_2}
- \\
& \
\ \ \ \ \ \ \ \ \ \
-
\frac{1}{4}\,\delta_{j_2}^{k_1}\,H_{j_1,j_1}^{j_1}\,\Theta^{n+1}
+
\frac{1}{4}\,\delta_{j_2}^{k_1}\,\sum_{k_2}\,
L_{j_1}^{k_2}\,\Theta^{k_2}
+
\frac{1}{4}\,\delta_{j_2}^{k_1}\,\Theta^{j_1}\,\Theta^{n+1}.
\endaligned\right.
\end{equation}
Sixthly:
\def\theequation{3.22}\begin{equation}
\left\{
\aligned
\delta_{j_1}^{k_1}\,\Theta_y^{n+1}
&
=
-
\underline{
L_{j_1, y}^{k_1}
+
2\,M_{x^{j_1}}^{k_1}
}
+ \\
& \
\ \ \ \ \ 
+
2\,\sum_{k_2}\,H_{j_1,k_2}^{k_1}\,M^{k_2}
-
\delta_{j_1}^{k_1}\,\sum_{k_2}\,H_{k_2,k_2}^{k_2}\,M^{k_2}
-
\frac{1}{2}\,\sum_{k_2}\,L_{j_1}^{k_2}\,L_{k_2}^{k_1}
+ \\
& \
\ \ \ \ \
+
\delta_{j_1}^{k_1}\,\sum_{k_2}\,M^{k_2}\,\Theta^{k_2}
+
\frac{1}{2}\,\delta_{j_1}^{k_1}\,\Theta^{n+1}\,\Theta^{n+1}.
\endaligned\right.
\end{equation}

\subsection*{3.23.~Solving $\Theta_{ x^{j_2 } }^{ j_1}$, 
$\Theta_y^{j_1}$, $\Theta_{x^{ j_1}}^{ n+1}$ and $\Theta_y^{ n+1 }$}
From the six families of equations~\thetag{ 3.17}, \thetag{ 3.18},
\thetag{ 3.19}, \thetag{ 3.20}, \thetag{ 3.21} and~\thetag{ 3.22}, we
can solve $\Theta_{ x^{j_2 }}^{j_1}$, $\Theta_y^{ j_1}$, $\Theta_{ x^{
j_1}}^{ n+1}$ and $\Theta_y^{ n+1}$. Not mentioning the (hard)
intermediate computations, we obtain firstly:
\def\theequation{3.24}\begin{equation}
\left\{
\aligned
\Theta_{x^{j_2}}^{j_1}
&
=
-
\underline{
2\,G_{j_1,j_2,y}
+
H_{j_1,j_1,x^{j_2}}^{j_1}
}
+
\sum_l\,G_{j_2,l}\,L_{j_1}^l
+
\frac{1}{2}\,H_{j_1,j_1}^{j_1}\,H_{j_2,j_2}^{j_2}
-
\sum_l\,H_{j_1,j_2}^l\,H_{l,l}^l
- \\
& \
\ \ \ \ \
-
G_{j_1,j_2}\,\Theta^{n+1}
-
\frac{1}{2}\,H_{j_1,j_1}\,\Theta^{j_1}
-
\frac{1}{2}\,H_{j_2,j_2}^{j_2}\,\Theta^{j_1}
+
\sum_l\,H_{j_1,j_2}^l\,\Theta^l
+
\frac{1}{2}\,\Theta^{j_1}\,\Theta^{j_2}.
\endaligned\right.
\end{equation}
Secondly:
\def\theequation{3.25}\begin{equation}
\left\{
\aligned
\Theta_y^{j_1}
&
=
-
\underline{
\frac{1}{3}\,H_{j_1,j_1,y}^{j_1}
+
\frac{2}{3}\,L_{j_1,x^{j_1}}^{j_1}
}
+
\frac{4}{3}\,G_{j_1,j_1}\,M^{j_1}
+
\frac{2}{3}\,\sum_l\,G_{j_1,l}\,M^l
-
\frac{1}{2}\,\sum_l\,H_{l,l}^l\,L_{j_1}^l
+ \\
& \
\ \ \ \ \
+
\frac{2}{3}\,\sum_l\,H_{j_1,l}^{j_1}\,L_{j_1}^l
-
\frac{2}{3}\,\sum_l\,H_{j_1,j_1}^l\,L_l^{j_1}
-
\frac{1}{2}\,H_{j_1,j_1}^{j_1}\,\Theta^{n+1}
+
\frac{1}{2}\,\sum_l\,L_{j_1}^l\,\Theta^l
+ \\
& \
\ \ \ \ \
+
\frac{1}{2}\,\Theta^{j_1}\,\Theta^{n+1}.
\endaligned\right.
\end{equation}
Thirdly:
\def\theequation{3.26}\begin{equation}
\left\{
\aligned
\Theta_{x^{j_1}}^{n+1}
&
=
-
\underline{
\frac{2}{3}\,H_{j_1,j_1,y}^{j_1}
+
\frac{1}{3}\,L_{j_1,x^{j_1}}^{j_1}
}
+
\frac{2}{3}\,G_{j_1,j_1}\,M^{j_1}
+
\frac{4}{3}\,\sum_l\,G_{j_1,l}\,M^l
-
\frac{1}{2}\,\sum_l\,H_{l,l}^l\,L_{j_1}^l
+ \\
& \
\ \ \ \ \ 
+
\frac{1}{3}\,\sum_l\,H_{j_1,l}^{j_1}\,L_{j_1}^l
-
\frac{1}{3}\,\sum_l\,H_{j_1,j_1}^l\,L_l^{j_1}
-
\frac{1}{2}\,H_{j_1,j_1}^{j_1}\,\Theta^{n+1}
+
\frac{1}{2}\,\sum_l\,L_{j_1}^l\,\Theta^l
+ \\
& \
\ \ \ \ \
+
\frac{1}{2}\,\Theta^{j_1}\,\Theta^{n+1}.
\endaligned\right.
\end{equation}
Fourtly:
\def\theequation{3.27}\begin{equation}
\left\{
\aligned
\Theta_y^{n+1}
&
=
-
\underline{
L_{j_1,y}^{j_1}
+
2\,M_{x^{j_1}}^{j_1}
}
+
2\,\sum_l\,H_{j_1,l}^{j_1}\,M^l
-
\sum_l\,H_{l,l}^l\,M^l
-
\frac{1}{2}\,\sum_l\,L_{j_1}^l\,L_l^{j_1}
+ \\
& \
\ \ \ \ \ 
+
\sum_l\,M^l\,\Theta^l
+
\frac{1}{2}\,\Theta^{n+1}\,\Theta^{n+1}.
\endaligned\right.
\end{equation}
These four families of partial differential equations constitute the
{\sl second auxiliary system}. By replacing these solutions in the
three remaining families of equations~\thetag{ 3.20}, \thetag{ 3.21}
and~\thetag{ 3.22}, we obtain supplementary equations
(which we do not copy) that are direct consequences of 
(I'), (II'), (III'), (IV').

To complete the proof of the main Lemma~3.3 above, it suffices now to
establish the first implication of the following list, since the other
three have been already established.

\begin{itemize}
\item[$\bullet$]
Some given functions $G_{j_1, j_2}$, $H_{j_1, j_2}^{k_1}$,
$L_{j_1}^{k_1}$ and $M^{k_1}$ of $(x^{l_1},y)$ satisfy the four
families of partial differential equations (I'), (II'), (III') and
(IV') of Theorem~1.1.
\item[$\Downarrow$]
\item[$\bullet$]
There exist functions $\Theta^{j_1}$, $\Theta^{n+1}$ satisfying the
second auxiliary system~\thetag{ 3.24}, \thetag{ 3.25}, \thetag{ 3.26}
and~\thetag{ 3.27}.
\item[$\Downarrow$]
\item[$\bullet$]
These solution functions $\Theta^{j_1}$, $\Theta^{n+1}$ satisfy the
six families of partial differential equations~\thetag{ 3.17},
\thetag{ 3.18}, \thetag{ 3.19}, \thetag{ 3.20}, \thetag{ 3.21}
and~\thetag{ 3.22}.
\item[$\Downarrow$]
\item[$\bullet$]
There exist functions $\Pi_{ j_1, j_2}^{k_1}$ of $(x^{l_1}, y)$, $1
\leq j_1, j_2, k_1 \leq m+1$, satisfying the first auxiliary
system~\thetag{ 3.7} of partial differential equations.
\item[$\Downarrow$]
\item[$\bullet$]
There exist functions $X^{l_2}$, $Y$ of $(x^{l_1},y)$ transforming the
system $y_{x^{j_1}x^{j_2}} = F^{j_1,j_2} (x^{l_1}, y, y_{x^{l_2}})$,
$j_1, j_2 = 1, \dots, n$, to the simplest system $Y_{X^{j_1}X^{j_2}} =
0$, $j_1, j_1 = 1, \dots, n$.
\end{itemize}

\smallskip

\subsection*{ 3.28.~Compatibility conditions for the second
auxiliary system} We notice that the second auxiliary system is also a
complete system. Thus, to establish the first above implication, it
suffices to show that the four families of compatibility conditions:
\def\theequation{3.29}\begin{equation}
\left\{
\aligned
0
& 
=
\left(
\Theta_{x^{j_2}}^{j_1}
\right)_{x^{j_3}}
-
\left(
\Theta_{x^{j_3}}^{j_1}
\right)_{x^{j_2}}, 
\\
0
& 
=
\left(
\Theta_{x^{j_2}}^{j_1}
\right)_y
-
\left(
\Theta_y^{j_1}
\right)_{x^{j_2}}, 
\\
0
& 
=
\left(
\Theta_{x^{j_1}}^{n+1}
\right)_{x^{j_2}}
-
\left(
\Theta_{x^{j_2}}^{n+1}
\right)_{x^{j_1}}, 
\\
0
& 
=
\left(
\Theta_{x^{j_2}}^{n+1}
\right)_y
-
\left(
\Theta_y^{n+1}
\right)_{x^{j_2}}, 
\endaligned\right.
\end{equation}
are a consequence of (I'), (I''), (III'), (IV').

For instance, in $(3.29)_1$, replacing $\Theta_{x^{j_2}}^{j_1}$ by its
expression~\thetag{ 3.24}, differentiating it with respect to
$x^{j_3}$, replacing $\Theta_{x^{j_3}}^{j_1}$ by its
expression~\thetag{ 3.24}, differentiating it with respect to
$x^{j_2}$ and substracting, we get: 
\def\theequation{3.30}\begin{equation}
\left\{
\aligned
0
& 
=
-
2\,G_{j_1,j_2,yx^{j_3}}
+
2\,G_{j_1,j_3,yx^{j_2}}
+
\underline{ 
H_{j_1,j_1,x^{j_2}x^{j_3}}^{j_1}
}_{ \octagon \! \! \! \! \! \tiny{\sf a}}
-
\underline{
H_{j_1,j_1,x^{j_3}x^{j_2}}^{j_1}
}_{ \octagon \! \! \! \! \! \tiny{\sf a}}
+ \\
& \
\ \ \ \ \
+
\frac{1}{2}\,
\underline{\Theta_{x^{j_3}}^{j_1}}\,\Theta^{j_2}
+
\frac{1}{2}\,\Theta^{j_1}\,
\underline{\Theta_{x^{j_3}}^{j_2}}
-
\frac{1}{2}\,
\underline{\Theta_{x^{j_2}}^{j_1}}\,\Theta^{j_3}
-
\frac{1}{2}\,\Theta^{j_1}\,
\underline{\Theta_{x^{j_2}}^{j_3}}
- \\
& \
\ \ \ \ \
-
\frac{1}{2}\,H_{j_1,j_1,x^{j_3}}^{j_1}\,\Theta^{j_2}
-
\frac{1}{2}\,H_{j_1,j_1}^{j_1}\,
\underline{\Theta_{x^{j_3}}^{j_2}}
+
\frac{1}{2}\,H_{j_1,j_1,x^{j_2}}^{j_1}\,\Theta^{j_3}
+
\frac{1}{2}\,H_{j_1,j_1}^{j_1}\,
\underline{\Theta_{x^{j_2}}^{j_3}}
- \\
& \
\ \ \ \ \ 
-
\frac{1}{2}\,H_{j_2,j_2,x^{j_3}}^{j_2}\,\Theta^{j_1}
-
\frac{1}{2}\,H_{j_2,j_2}^{j_2}\,
\underline{\Theta_{x^{j_3}}^{j_1}}
+
\frac{1}{2}\,H_{j_3,j_3,x^{j_2}}^{j_3}\,\Theta^{j_1}
+
\frac{1}{2}\,H_{j_3,j_3}^{j_3}\,
\underline{\Theta_{x^{j_2}}^{j_1}}
- \\
& \
\ \ \ \ \
-
G_{j_1,j_2,x^{j_3}}\,\Theta^{n+1}
-
G_{j_1,j_2}\,\underline{\Theta_{x^{j_3}}^{n+1}}
+
G_{j_1,j_3,x^{j_2}}\,\Theta^{n+1}
+
G_{j_1,j_3}\,\underline{\Theta_{x^{j_2}}^{n+1}}
+ \\
& \
\ \ \ \ \
+
\sum_l\,H_{j_1,j_2,x^{j_3}}^l\,\Theta^l
+
\sum_l\,H_{j_1,j_2}^l\,
\underline{\Theta_{x^{j_3}}^l}
-
\sum_l\,H_{j_1,j_3,x^{j_2}}^l\,\Theta^l
-
\sum_l\,H_{j_1,j_3}^l\,
\underline{\Theta_{x^{j_2}}^l}
+ \\
& \
\ \ \ \ \
+
\frac{1}{2}\,H_{j_1,j_1,x^{j_3}}^{j_1}\,H_{j_2,j_2}^{j_2}
+
\frac{1}{2}\,H_{j_1,j_1}^{j_1}\,H_{j_2,j_2,x^{j_3}}^{j_2}
-
\frac{1}{2}\,H_{j_1,j_1,x^{j_2}}^{j_1}\,H_{j_3,j_3}^{j_3}
-
\frac{1}{2}\,H_{j_1,j_1}^{j_1}\,H_{j_3,j_3,x^{j_2}}^{j_3}
- \\
& \
\ \ \ \ \ 
-
\sum_l\,H_{j_1,j_2,x^{j_3}}^l\,H_{l,l}^l
-
\sum_l\,H_{j_1,j_2}^l\,H_{l,l,x^{j_3}}^l
+
\sum_l\,H_{j_1,j_3,x^{j_2}}^l\,H_{l,l}^l
+
\sum_l\,H_{j_1,j_3}^l\,H_{l,l,x^{j_2}}^l
+ \\
& \
\ \ \ \ \ 
+
\sum_l\,G_{j_2,l,x^{j_3}}\,L_{j_1}^l
+
\sum_l\,G_{j_2,l}\,L_{j_1,x^{j_3}}^l
-
\sum_l\,G_{j_3,l,x^{j_2}}\,L_{j_1}^l
-
\sum_l\,G_{j_3,l}\,L_{j_1,x^{j_2}}^l.
\endaligned\right.
\end{equation}
Next, replacing the twelve first order partial derivatives underlined
just above:
\def\theequation{3.31}\begin{equation}
\left\{
\aligned
&
\underline{\Theta_{x^{j_3}}^{j_1}}, 
\ \ \ \ \
\underline{\Theta_{x^{j_3}}^{j_2}}, 
\ \ \ \ \
\underline{\Theta_{x^{j_2}}^{j_1}}, 
\ \ \ \ \
\underline{\Theta_{x^{j_2}}^{j_3}},
\ \ \ \ \
\underline{\Theta_{x^{j_3}}^{j_2}}, 
\ \ \ \ \
\underline{\Theta_{x^{j_2}}^{j_3}},
\\
&
\underline{\Theta_{x^{j_3}}^{j_1}},
\ \ \ \ \
\underline{\Theta_{x^{j_2}}^{j_1}},
\ \ \ \ \
\underline{\Theta_{x^{j_3}}^{n+1}},
\ \ \ \ \
\underline{\Theta_{x^{j_2}}^{n+1}},
\ \ \ \ \
\underline{\Theta_{x^{j_3}}^l},
\ \ \ \ \
\underline{\Theta_{x^{j_2}}^l}.
\endaligned\right.
\end{equation}
by their values issued from~\thetag{ 3.24}, \thetag{ 3.26} and
adapting the summation indices, we get the explicit developed form of
the first family of compatibility conditions $(3.29)_1$:
\def\theequation{3.32}\begin{equation}
\left\{
\aligned
0 
= 
?
&
=
-
\underline{\underline{
2\,G_{j_1,j_2,x^{j_3}y}
+
2\,G_{j_1,j_3,x^{j_2}y}
}}
- \\
& \
\ \ \ \ \
-
\underline{
\sum_l\,G_{j_3,l,x^{j_2}}\,L_{j_1}^l
+
\sum_l\,G_{j_2,l,x^{j_3}}\,L_{j_1}^l
-
G_{j_1,j_2,y}\,H_{j_3,j_3}^{j_3}
+
G_{j_1,j_3,y}\,H_{j_2,j_2}^{j_2}
}
- \\
& \
\ \ \ \ \
-
\underline{
2\,\sum_l\,G_{l,j_3}\,H_{j_1,j_2}^l
+
2\,\sum_l\,G_{l,j_2}\,H_{j_1,j_3}^l
-
\sum_l\,H_{j_1,j_2,x^{j_3}}^l\,H_{l,l}^l
+
\sum_l\,H_{j_1,j_3,x^{j_2}}^l\,H_{l,l}^l
}
- \\
& \
\ \ \ \ \
-
\underline{
\frac{2}{3}\,H_{j_2,j_2,y}^{j_2}\,G_{j_1,j_3}
+
\frac{2}{3}\,H_{j_3,j_3,y}^{j_3}\,G_{j_1,j_2}
-
\frac{2}{3}\,L_{j_3,x^{j_3}}^{j_3}\,G_{j_1,j_2}
+
\frac{2}{3}\,L_{j_2,x^{j_2}}^{j_2}\,G_{j_1,j_3}
}
- \\
& \
\ \ \ \ \
-
\underline{
\sum_l\,L_{j_1,x^{j_2}}^l\,G_{j_3,l}
+
\sum_l\,L_{j_1,x^{j_3}}^l\,G_{j_2,l}
}
- \\
& \
\ \ \ \ \
-
\frac{2}{3}\,G_{j_1,j_2}\,G_{j_3,j_3}\,M^{j_3}
+
\frac{2}{3}\,G_{j_1,j_3}\,G_{j_2,j_2}\,M^{j_2}
-
\frac{4}{3}\,\sum_l\,G_{j_1,j_2}\,G_{j_3,l}\,M^l
+ \\
& \
\ \ \ \ \
+
\frac{4}{3}\,\sum_l\,G_{j_1,j_3}\,G_{j_2,l}\,M^l
-
\frac{1}{2}\,\sum_l\,G_{j_3,l}\,H_{j_1,j_1}^{j_1}\,L_{j_2}^l
+
\frac{1}{2}\,\sum_l\,G_{j_2,l}\,H_{j_1,j_1}^{j_1}\,L_{j_3}^l
- \\
& \
\ \ \ \ \
-
\frac{1}{2}\,\sum_l\,G_{j_3,l}\,H_{j_2,j_2}^{j_2}\,L_{j_1}^l
+
\frac{1}{2}\,\sum_l\,G_{j_2,l}\,H_{j_3,j_3}^{j_3}\,L_{j_1}^l
-
\frac{1}{2}\,\sum_l\,G_{j_1,j_3}\,H_{l,l}^l\,L_{j_2}^l
+ \\
& \
\ \ \ \ \
+
\frac{1}{2}\,\sum_l\,G_{j_1,j_2}\,H_{l,l}^l\,L_{j_3}^l
-
\frac{1}{3}\,\sum_l\,G_{j_1,j_2}\,H_{j_3,l}^{j_3}\,L_{j_3}^l
+
\frac{1}{3}\,\sum_l\,G_{j_1,j_3}\,H_{j_2,l}^{j_2}\,L_{j_2}^l
- \\
& \
\ \ \ \ \
-
\frac{1}{3}\,G_{j_1,j_3}\,H_{j_2,j_2}^l\,L_l^{j_2}
+
\frac{1}{3}\,G_{j_1,j_2}\,H_{j_3,j_3}^l\,L_l^{j_3}
- \\
& \
\ \ \ \ \
-
\sum_l\,\sum_p\,G_{j_2,p}\,H_{j_1,j_3}^l\,L_l^p
+
\sum_l\,\sum_p\,G_{j_3,p}\,H_{j_1,j_2}^l\,L_l^p
- \\
& \
\ \ \ \ \
-
\sum_l\,\sum_p\,H_{j_1,j_2}^l\,H_{l,j_3}^p\,H_{p,p}^p
+
\sum_l\,\sum_p\,H_{j_1,j_3}^l\,H_{l,j_2}^p\,H_{p,p}^p.
\endaligned\right.
\end{equation}

\def\thelemma{3.33}\begin{lemma}
{\rm (\cite{ m2003}, \cite{ m2004a})}
This first family of compatibility conditions for the second auxiliary
system obtained by developing $(3.29)_1$ in length, together with the
three remaining families obtained by developing $(3.29)_2$,
$(3.29)_3$, $(3.29)_4$ in length, are consequences, by linear
combinations and by differentiations, of {\rm (I')}, {\rm (II')}, {\rm
(III')}, {\rm (IV')}, of Theorem~1.1.
\end{lemma}

The proof of Theorem~1.1 is complete.
\qed

\vfill
\begin{thebibliography}{GTW1989}

\bibitem[Bi2003]{bi2003}
{\sc Bi\`eche}, C.:
{\em Le probl\`eme d'\'equivalence locale pour un syst\`eme scalaire
complet d'\'equations aux d\'eriv\'ees partielles d'ordre deux \`a $n$
variables ind\'ependantes}, Preprint, 25~pp., July 2003.

\bibitem[BN2002]{bn2002}
{\sc Bi\`eche}, C.; {\sc Neut}, S.:
{\em Sur le probl\`eme d'\'equivalence de certains syst\`emes
d'\'equations aux d\'eriv\'ees partielles d'ordre 2}, 
Preprint, 2002.

\bibitem[BK1989]{bk1989}
{\sc Bluman}, G.W.; {\sc Kumei}, S.:
{\em Symmetries and differential equations}, 
Springer-Verlag, New-York, 1989.

\bibitem[Ca1924]{ca1924}
{\sc Cartan}, \'E.:
{\em Sur les vari\'et\'es \`a connexion projective}, 
Bull. Soc. Math. France {\bf 52} (1924), 205--241.

\bibitem[CM1974]{cm1974}
{\sc Chern}, S.S.; {\sc Moser}, J.K.:
{\em Real hypersurfaces in complex manifolds},
Acta Math. {\bf 133} (1974), no.~2, 219--271.

\bibitem[Ch1975]{ch1975}
{\sc Chern}, S.-S.:
{\em On the projective structure of a real hypersurface in 
$\C^{n+1}$}, Math. Scand. {\bf 36} (1975), 74--82.

\bibitem[Fe1995]{fe1995}
{\sc Fels}, M.:
{\em The equivalence problem for systems of 
second-order ordinary differential equations}, Proc. 
London Math. Soc. {\bf 71} 
(1995), no.~2, 221--240.

\bibitem[G1989]{g1989}
{\sc Gardner}, R.B.:
{\em The method of equivalence and its applications}, 
CBMS-NSF Regional Conference Series in Applied Mathematics 
{\bf 58} (SIAM, Philadelphia, 1989), 127~pp.

\bibitem[GTW1989]{gtw1989}
{\sc Grissom}, C.; {\sc Thompson}, G.; {\sc Wilkens}, G.:
{\em Linearization of second order ordinary differential equations 
via Cartan's equivalence method}, J. Diff. Eq. {\bf 77}
(1989), no.~1, 1--15.

\bibitem[GM2003]{gm2003}
{\sc Gaussier}, H.; {\sc Merker}, J.: 
{\em Symmetries of partial differential equations},
J. Korean Math. Soc. {\bf 40} (2003), no.~3, 517--561.

\bibitem[Ha1937]{ha1937}
{\sc Hachtroudi}, M.:
{\em Les espaces d'\'el\'ements \`a connexion projective normale},
Actualit\'es Scientifiques et Industrielles, 565, Paris, Hermann, 
1937. 

\bibitem[HK1989]{hk1989}
{\sc Hsu}, L.; {\sc Kamran}, N.:
{\em Classification of second order ordinary differential equations
admitting Lie groups of fibre-preserving point symmetries},
Proc. London Math. Soc. {\bf 58} (1989), no.~3, 387--416.

\bibitem[Lie1883]{lie1883}
{\sc Lie}, S.: 
{\em Klassifikation und Integration vo gew\"ohnlichen
Differentialgleichungen zwischen $x$, $y$, die eine
Gruppe von Transformationen gestaten I-IV}. In:
Gesammelte Abhandlungen, Vol.~5, B.G. Teubner, Leipzig, 1924, 
pp.~240--310; 362--427, 432--448.

\bibitem[M2003]{m2003}
{\sc Merker}, J.:
{\em hand manuscript I}, 212~pp., autumn 2003;
{\em hand manuscript II}, 114~pp., autumn 2003.

\bibitem[M2004a]{m2004a}
{\sc Merker}, J.:
{\em Explicit differential characterization of the newtonian free
particle system in $m\geq 2$ dependent variables}, Preprint,
University of Provence, November 2004.

\bibitem[M2004b]{m2004b}
{\sc Merker}, J.:
{\em Four explicit formulas for the $\kappa$-th point prolongation 
of a Lie symmetry in an arbitrary number of independent and of dependent
variables}, 
in preparation.

\bibitem[M2004c]{m2004c}
{\sc Merker}, J.:
{\em Symmetries of completely integrable systems of 
analytic partial differential equations},
in preparation.

\bibitem[M2005]{m2005}
{\sc Merker}, J.:
{\em Explicit Cartan-Hachtroudi-Chern tensors}, 
in preparation.

\bibitem[N2003]{n2003}
{\sc Neut}, S.:
{\em Implantation et nouvelles applications de la m\'ethode
d'\'equivalence d'\'Elie Cartan}, Th\`ese, Universit\'e 
Lille~1, October 2003.

\bibitem[OL1986]{ol1986} 
{\sc Olver},~P.J.:
{\em Applications of Lie groups to differential equations}.
Springer Verlag, Heidelberg, 1986.

\bibitem[OL1995]{ol1995} 
{\sc Olver},~P.J.:
{\em Equivalence, Invariance and Symmetries}. Cambridge
University Press, Cambridge, 1995, xvi+525~pp.

\bibitem[Sh1997]{sh1997} 
{\sc Sharpe},~R.W.:
{\em Differential geometry; Cartan's generalization of Klein's
Erlangen program}, 
Graduate texts in mathematics, vol.~166, 
Springer Verlag, Berlin, 1997, xix+421~pp.

\bibitem[Su2001]{su2001} 
{\sc Sukhov},~A.:
{\em Segre varieties and Lie symmetries},
Math. Z. {\bf 238} (2001), no.~3, 483--492.

\end{thebibliography}
\end{document}